# TENSOR PRODUCT REPRESENTATIONS FOR ORTHOSYMPLECTIC LIE SUPERALGEBRAS

**Georgia Benkart***
Department of Mathematics
University of Wisconsin
Madison, WI 53706 USA

**Chanyoung Lee Shader****
Department of Mathematics
University of Wyoming
Laramie, WY 82071 USA

**Arun Ram*****
Department of Mathematics
University of Sydney
NSW 2006 Australia

Abstract

We derive a general result about commuting actions on certain objects in braided rigid monoidal categories. This enables us to define an action of the Brauer algebra on the tensor space $V^{\otimes k}$ which commutes with the action of the orthosymplectic Lie superalgebra $\mathfrak{spo}(V)$ and the orthosymplectic Lie color algebra $\mathfrak{spo}(V,\beta)$. We use the Brauer algebra action to compute maximal vectors in $V^{\otimes k}$ and to decompose $V^{\otimes k}$ into a direct sum of submodules $T^\lambda$. We compute the characters of the modules $T^\lambda$, give a combinatorial description of these characters in terms of tableaux, and model the decomposition of $V^{\otimes k}$ into the submodules $T^\lambda$ with a Robinson-Schensted-Knuth type insertion scheme.

Outline

---

* Supported in part by National Science Foundation Grant DMS-9300523
** Supported in part by a University of Wyoming Basic Research Grant
*** Supported in part by National Science Foundation Grant DMS-9300523, by a National Science Foundation Postdoctoral Fellowship, and by an Australian Research Council Research Fellowship.

## 0. Introduction

*Summary of results*

In this paper we show that there are orthosymplectic Lie superalgebra and orthosymplectic Lie color algebra analogues of the results developed by Berele and Regev [BR] and Sergeev [Se] for general linear Lie superalgebras. Our work corresponds in the superalgebra and color algebra setting to what Brauer did in [Br] by extending to orthogonal and symplectic groups Schur's classical results for general linear groups.

In his thesis [Sh1] and a subsequent paper [Sh2] Schur proved that the action of the symmetric group $S_k$ on tensor space $V^{\otimes k}$ by place permutations determines the centralizer of the action of the general linear group $GL(V)$ on $V^{\otimes k}$. This result, often called *Schur-Weyl duality*, relates in a very fundamental way the representation theory and combinatorics of the groups $S_k$ and $GL(V)$. The orthogonal group $O(V)$ also acts on $V^{\otimes k}$, and Brauer [Br] constructed an algebra, now referred to as the *Brauer algebra*, which commutes with the $O(V)$-action on tensor space. When $V$ is even-dimensional, the symplectic group $Sp(V)$ also has an action on $V^{\otimes k}$, and its centralizer can be described using the Brauer algebra.

Berele and Regev [BR] and Sergeev [Se] have shown that the action of the symmetric group $S_k$ on $V^{\otimes k}$ by "graded" place permutations determines the centralizer of the general linear Lie superalgebra when $V$ is $\mathbb{Z}_2$-graded, and they have exploited this action to study certain modules for the superalgebra and their characters. Fischman and Montgomery [FM] have generalized the work of [BR] and [Se] to cotriangular Hopf algebras which arise from enveloping algebras of general linear Lie color algebras.

Our intent is to establish an orthosymplectic version of this theory. In particular, we

(i) show there is a Brauer algebra action on $V^{\otimes k}$ which commutes with the orthosymplectic Lie color algebra $\mathfrak{spo}(V,\beta)$, (the orthosymplectic Lie superalgebra $\mathfrak{spo}(V)$-action is just a special case);

(ii) use the Brauer algebra to construct a family of maximal vectors for $\mathfrak{spo}(V,\beta)$ in the tensor space (Theorem 3.9);

(iii) obtain from the Brauer algebra action a direct sum decomposition of $V^{\otimes k}$ into $\mathfrak{spo}(V,\beta)$-submodules $T^\lambda$ naturally indexed by partitions $\lambda$;

(iv) use the character theory of the Brauer algebra to compute the characters of the $\mathfrak{spo}(V,\beta)$-modules $T^\lambda$ and to give a combinatorial description of these characters in terms of tableaux;

(v) develop an insertion scheme for the tableaux which models the decomposition of the tensor space into the modules $T^\lambda$.

*Remarks on the results in this paper*

(a) To our knowledge, orthosymplectic Lie color algebras were first introduced in [B], which discusses the Brauer algebra action, but does not prove that this action commutes with $\mathfrak{spo}(V,\beta)$. The notion of an orthosymplectic Lie color algebra allows us to give a uniform proof that the Brauer algebra action on tensor space commutes with the action of the orthogonal group (Lie algebra), the symplectic group (Lie algebra), the orthosymplectic Lie superalgebra as well as more general group graded orthosymplectic algebras. As we discuss in Section 1, orthosymplectic Lie color algebras have a natural root space decomposition and a triangular decomposition which are exactly analogous to the Lie superalgebra case. The derivation of the maximal vectors which we give in Section 3 extends the work in [BBL] where maximal vectors are computed for classical case of the orthogonal and symplectic Lie algebras.

(b) The action of the Brauer algebra arises naturally from the structure of the category of finite dimensional modules for Lie color algebras. As we show in Section 2, braided monoidal categories provide a convenient framework for developing the results on commuting actions. Our arguments in Section 2 will apply to give a centralizing action of an algebra (in general it may be not be a Brauer algebra) on tensor space for any kind of "Lie like" algebra or quantum group for which the category of finite dimensional modules has a braided monoidal structure and a special isomorphism between $V$ and $V^*$. Thus for example, one can modify our Section 2 with very minor changes to obtain an action of the Birman-Wenzl-Murakami algebra centralizing the appropriate quantum analogue of $\mathfrak{spo}(V,\beta)$. The idea is essentially the same as that used by Reshetikhin-Turaev [RT], in particular, see for example [CP, Theorem 5.3.8]. It shares similarities with the methods that Fischman and Montgomery [FM] apply to get their symmetric group action, except that we have chosen to work with the structure on the category of modules rather than with the Hopf algebra structure as they do.

(c) Berele and Regev [BR] give a very interesting combinatorial description of the characters of the $\mathfrak{gl}(m/n)$-modules which appear in the tensor space $V^{\otimes k}$ ($V = V_0 \oplus V_1$, $\dim V_0 = m$, and $\dim V_1 = n$) by describing them as hybrid Schur functions involving tableaux which have a column-strict part and a row-strict part. In a similar fashion we give a combinatorial description of the character of the module $T^\lambda$ as a hybrid symplectic-ordinary Schur function given by tableaux which have a symplectic part and a row-strict part. The symplectic part is a symplectic tableau of the kind introduced by King (see for example, [KE]) to index basis elements for irreducible representations of symplectic groups.

(d) The modules $T^\lambda$ are the same as the ones considered by Bars and Balantekin in [BB1-BB2]. They give Jacobi-Trudi type character formulas, but they do not derive the tableaux description that we give here. Bars and Balentekin seem to indicate that the modules $T^\lambda$ are irreducible, but this is not clear to us, either from their work or from ours. In fact, R.C. King in personal communication has told us that he has found explicit examples of $T^\lambda$ which are not irreducible.

(e) There has been other work, notably ([FJ], [CK], [L1-L2]), which describes how to index representations of the Lie superalgebra $\mathfrak{spo}(V)$ by partitions, but none of these papers has given an interpretation for their characters in terms of tableaux. The main ingredient in developing the tableaux description is the identity in Theorem 4.24 (h). A very similar version of this identity (Theorem 4.24 (i)) appears in [CK]. This identity could be used in combination with the work of Sundaram [Su2] to give another combinatorial interpretation for these characters. See [Ki2] for a survey of the use of tableaux in the study of representations of Lie superalgebras.

(f) There is an extensive literature of papers by Bernstein and Leites [BL], [Le1-Le2], Kac [K1-K2], van der Jeugt, Hughes, King, and Thierry-Mieg [JHKT1-JHKT3], Penkov and Serganova [PS1-PS3], [P], Kac and Wakimoto [KW], and others which studies representations of Lie superalgebras using Lie theoretic and geometric methods. These approaches also yield character formulas, the most general of which is the Weyl-Kac character formula. We have not made any effort to understand our character formulas in this other setting, although the formulas must be equal in many cases. Even for the superalgebra $\mathfrak{gl}(m/n)$, the connection between the results of [BR] and [Se], the Sergeev-Pragacz character formula, and the Weyl-Kac character formula needs to be better understood. King and others have done some work in this direction, (see [Ki2]). The relationship between the centralizer approach to the representation theory of Lie superalgebras and the approach via Kac modules and typical weights also needs to be better explained.

(g) In our work we have not included proofs of the analogue of Schur-Weyl duality, i.e. we have not shown that the actions of the Brauer algebra and the Lie color algebra $\mathfrak{spo}(V, \beta)$ each generate the full centralizer of the action of the other. Let us make just a few remarks on this point.

(i) In the general linear-symmetric group case of [BR], [Se], and [FM], one can get away with proving only half of the duality and using the semisimplicity of the group algebra of the symmetric group $\mathbb{C}S_k$ to obtain the other half for free. This is not possible in the orthosymplectic-Brauer algebra case since the Brauer algebra is not necessarily semisimple.

(ii) The usual trick for proving that general linear lie algebra $\mathfrak{gl}(V)$ generates the full centralizer of the symmetric group action uses the idempotent $\sum_{\sigma\in S_k}\sigma$ (which corresponds to the trivial $S_k$-character) to construct a projection map onto the $\mathfrak{gl}(V)$-invariants. Unfortunately, this is not available in the Brauer algebra case since the Brauer algebra does not have a one-dimensional module which affords invariants.

(iii) We have succeeded in establishing various parts of both halves of the duality in our orthosymplectic-Brauer algebra setting but have chosen not to include these results in this paper. Optimally what one would like is a proof of the duality which handles all the cases simultaneously rather than arguing separately for the orthogonal group, then the symplectic group, then the orthosymplectic Lie superalgebras, then the quantum orthogonal group, etc.

*Open problems*

(1) The relationship between the centralizer approach to the representation theory of Lie superalgebras and the approach using Kac modules and typical-atypical weights needs to be better understood. To our knowledge this is only partially done even for the $\mathfrak{gl}(m/n)$ case [JHKT1-3], and it is not known how to verify directly that the Weyl-Kac character formula for typical representations and the Pragacz-Sergeev character formula for $\mathfrak{gl}(m/n)$-irreducible modules in tensor space are equal. Here we have not made any attempt to relate our results to the Weyl-Kac formula, although this should be done sometime in the future.

(2) In determining the characters of the modules $T^\lambda$, we have shown that they are equal to polynomials $sc_\mu(x_1,x_1^{-1},\dots,x_r,x_r^{-1},y_1,y_1^{-1},\dots,y_s,y_s^{-1},1)$ which appear as the coefficients of the Schur functions $s_\mu(z_1,\dots,z_{r+s})$ in the identity

$$\prod_{1\le i<j\le r+s}(1-z_iz_j)\frac{\prod_{j=1}^{r+s}(1+z_j)\prod_{j=1}^{r+s}\prod_{i=1}^{s}(1+y_iz_j)(1+y_i^{-1}z_j)}{\prod_{j=1}^{r+s}\prod_{i=1}^{r}(1-x_iz_j)(1-x_i^{-1}z_j)}$$
$$=\sum_\mu sc_\mu(x_1,x_1^{-1},\dots,x_r,x_r^{-1},y_1,y_1^{-1},\dots,y_s,y_s^{-1},1)s_\mu(z_1,\dots,z_{r+s}). \tag{0.1}$$

There are two classical identities of Littlewood [Li] and Weyl [We] for the characters $sp_\mu$ of the symplectic group $Sp(2r)$ and the characters $so_\mu$ of the orthogonal groups $SO(2s+1)$:

$$\prod_{1\le i<j\le r}(1-z_iz_j)\frac{1}{\prod_{j=1}^{r}\prod_{i=1}^{r}(1-x_iz_j)(1-x_i^{-1}z_j)}$$
$$=\sum_\mu sp_\mu(x_1,x_1^{-1},\dots,x_r,x_r^{-1})s_\mu(z_1,\dots,z_r), \tag{0.2}$$

$$\prod_{1\le i<j\le s}(1-z_iz_j)\prod_{j=1}^{s}(1+z_j)\prod_{j=1}^{s}\prod_{i=1}^{s}(1+y_iz_j)(1+y_i^{-1}z_j)$$
$$=\sum_\mu so_{\mu'}(y_1,y_1^{-1},\dots,y_s,y_s^{-1},1)s_\mu(z_1,\dots,z_s), \tag{0.3}$$

where $\mu'$ is the conjugate of the partition $\mu$. When the orthogonal portion is zero, i.e. $s = 0$ and the variable 1 is not present, identity (0.1) gives the classical result in (0.2), and when the symplectic part is zero, i.e. $r = 0$, identity (0.1) reduces to (0.3). Is there a combinatorial interpretation for the functions $sc_\mu$ which expresses $sc_\mu$ as a hybrid object built from symplectic and orthogonal characters? Our combinatorial description of $sc_\mu$ is as a mixed symplectic-general linear character rather than as a hybrid symplectic-orthogonal character.

(3) Find a general Schur-Weyl duality result, see remark (g) above.

*Acknowledgements*

The investigations in this paper began with a study of the maximal vectors and stability behaviour of $V^{\otimes k}$ for the orthosymplectic Lie superalgebra $\mathfrak{spo}(V)$ in [L1]. The other results in this paper were obtained primarily in 1994 and early 1995 and the paper was finalized in late 1995. All of us are very grateful for the support our research has received from various grants and institutions during this time. We thank the National Science Foundation, the University of Wyoming, and the Australian Research Council for support at various stages of the production of this paper, as follows:

Georgia Benkart: National Science Foundation Grant DMS-9300523.

Chanyoung Lee Shader: University of Wyoming Basic Research Grant.

Arun Ram: National Science Foundation grant DMS-9300523, an NSF Postdoctoral Fellowship and an Australian Research Council Research Fellowship.

## 1. Lie color algebras and $\mathfrak{spo}(V, \beta)$

*Lie color algebras*

Let $\kappa$ denote a field of characteristic zero (this assumption could be relaxed, but for convenience we stay in characteristic 0). Let $G$ be a finite abelian group with identity $1_G$. A *symmetric bicharacter* on $G$ is a map $\beta: G \times G \to \kappa^*$ into the multiplicative group of the field such that

(1) $\beta(ab, c) = \beta(a, c)\beta(b, c)$,

(2) $\beta(a, bc) = \beta(a, b)\beta(a, c)$, and

(3) $\beta(a, b)\beta(b, a) = 1$, for all $a, b \in G$.

Taking $a = 1_G$ in the first relation and $b = 1_G$ in the second shows that $\beta(1_G, c) = 1 = \beta(a, 1_G)$ for all $a, c \in G$.

A $\kappa$-vector space $V$ is *$G$-graded* if it is the direct sum $V = \bigoplus_{a \in G} V_a$ of subspaces indexed by the elements of $G$. If $v \in V_a$ for some $a \in G$, then $v$ is *homogeneous of degree $a$*.

Assume $\beta$ is a fixed symmetric bicharacter on a group $G$. A *Lie color algebra* $(\mathfrak{g}, G, \beta)$ is a $G$-graded vector space $\mathfrak{g} = \bigoplus_{a \in G} \mathfrak{g}_a$ with a $\kappa$-bilinear bracket $[\,,\,]: \mathfrak{g} \times \mathfrak{g} \to \mathfrak{g}$ such that

(1) $[\mathfrak{g}_a, \mathfrak{g}_b] \subseteq \mathfrak{g}_{ab}$, for all $a, b \in G$,

(2) $[x, y] = -\beta(b, a)[y, x]$, for $x \in \mathfrak{g}_a$, $y \in \mathfrak{g}_b$, and

(3) $[x, [y, z]] = [[x, y], z] + \beta(b, a)[y, [x, z]]$, for $x \in \mathfrak{g}_a$, $y \in \mathfrak{g}_b$, and all $z \in \mathfrak{g}$.

When the group is the cyclic group $G = \{\pm 1\} = \{(-1)^a \mid a = 0,1\}$ order 2 and $\beta((-1)^a,(-1)^b) = (-1)^{ab}$, then $\mathfrak{g}$ is a *Lie superalgebra*. In this case it is customary to regard the grading as by the additive group of the field $\mathbb{Z}_2 = \{0,1\}$, $\mathfrak{g} = \mathfrak{g}_0 \oplus \mathfrak{g}_1$, with $\beta(a,b) = (-1)^{ab}$ for $a,b \in \mathbb{Z}_2$.

*Remark 1.1.* Since, for a Lie color algebra, the bicharacter $\beta$ is symmetric, $\beta(a^{-1},a) = \beta(a,a)^{-1} = \beta(a,a) = \pm 1$. Thus $\mathfrak{g} = \mathfrak{g}_{(0)} \oplus \mathfrak{g}_{(1)}$, where

$$\mathfrak{g}_{(0)} = \bigoplus_{a,\,\beta(a,a)=1} \mathfrak{g}_a \quad \text{and} \quad \mathfrak{g}_{(1)} = \bigoplus_{a,\,\beta(a,a)=-1} \mathfrak{g}_a.$$

One shows easily that $\beta(ab,ab) = \beta(a,a)\beta(b,b)$, which implies the map $a \mapsto \beta(a,a) \in \{\pm 1\}$ is a homomorphism. It follows that the multiplication in $\mathfrak{g} = \mathfrak{g}_{(0)} \oplus \mathfrak{g}_{(1)}$ is $\mathbb{Z}_2$-graded, but still $\mathfrak{g}$ may *not* be a Lie superalgebra.

*The category of finite dimensional modules for a Lie color algebra*

Let $\mathfrak{g}$ be a Lie color algebra. A $\mathfrak{g}$-*module* is $G$-graded vector space $V = \bigoplus_{a\in G} V_a$ with a $\mathfrak{g}$-action $\mathfrak{g} \otimes V \to V$, $x \otimes v \mapsto xv$ such that

(1) if $x \in \mathfrak{g}_a$ and $v \in V_b$, then $xv \in V_{ab}$,

(2) $[x,y]v = x(yv) - \beta(b,a)y(xv)$, for all $x \in \mathfrak{g}_a, y \in \mathfrak{g}_b$ and all $v \in V$.

A $\mathfrak{g}$-*module morphism* from a $\mathfrak{g}$-module $M$ to a $\mathfrak{g}$-module $N$ is a $\kappa$-linear map $\phi\colon M \to N$ such that

(1) $\phi(M_a) \subseteq N_a$ for all $a \in G$, and

(2) $\phi(xm) = x\phi(m)$ for all $x \in \mathfrak{g}$ and all $m \in M$.

If $M$ and $N$ are $\mathfrak{g}$-modules, then the tensor product $M \otimes N$ is a $\mathfrak{g}$-module with

(1) $G$-grading given by $(M \otimes N)_c = \bigoplus_{ab=c} M_a \otimes N_b$, and

(2) $\mathfrak{g}$-action defined by

$$x(m \otimes n) = xm \otimes n + \beta(b,a)(m \otimes xn), \tag{1.2}$$

for all $x \in \mathfrak{g}_a$, $m \in M_b$, and $n \in N$.

If $f\colon M \to U$ and $g\colon N \to V$ are $\mathfrak{g}$-module morphisms, then

$$\begin{array}{rccc} f \otimes g\colon & M \otimes N & \longrightarrow & U \otimes V \\ & m \otimes n & \longmapsto & f(m) \otimes g(n) \end{array}$$

is a $\mathfrak{g}$-module morphism.

The *braiding morphisms* are the $\mathfrak{g}$-module isomorphisms defined by

$$\begin{array}{rccc} \check{R}_{M,N}\colon & M \otimes N & \longrightarrow & N \otimes M \\ & m \otimes n & \longmapsto & \beta(b,a) n \otimes m \end{array} \tag{1.3}$$

for all $m \in M_a$, $n \in N_b$.

The *trivial* $\mathfrak{g}$-*module* is the one-dimensional $\kappa$-vector space $\mathbf{1} = \kappa$ with grading $\mathbf{1} = \mathbf{1}_{1_G}$ and $\mathfrak{g}$-action defined by $xv = 0$ for all $x \in \mathfrak{g}$ and all $v \in \mathbf{1}$. The *dual* of a $\mathfrak{g}$-module $M$ is the vector space $M^* = \mathrm{Hom}_\kappa(M,\kappa)$ with

(1) $G$-grading given by $(M^*)_a = \{f \in M^* \mid f(M_b) = 0 \text{ if } b \neq a^{-1}\}$, and

(2) $\mathfrak{g}$-action given by

$$(xf)(m) = -\beta(b,a)f(xm), \tag{1.4}$$

for all $x \in \mathfrak{g}_a$, $f \in (M^*)_b$, and $m \in M$.

It follows easily by direct verification that the category of finite dimensional $\mathfrak{g}$-modules with the above constructions satisfies the definitions and axioms of a braided rigid monoidal category (see [CP] §5.2B). Although we phrase our results in this section and the next in the language of braided monoidal categories, the reader need not know anything about braided monoidal categories beyond the above standard constructions of modules. To summarize, we have

**Proposition 1.5.** *Let $(\mathfrak{g}, G, \beta)$ be a Lie color algebra. Then with the above constructions the category of finite dimensional $\mathfrak{g}$-modules is a braided strict rigid monoidal category.*

*The general linear Lie color algebra $\mathfrak{gl}(V,\beta)$*

Let $G$ be a finite abelian group with symmetric bicharacter $\beta$. Assume $V = \bigoplus_{a\in G} V_a$ is a $G$-graded $\kappa$-vector space. Let $\mathfrak{gl}(V,\beta) = \mathrm{End}(V)$ denote the $\kappa$-vector space of $\kappa$-linear maps from $V$ to $V$ with the $G$-grading assigned by

$$\mathfrak{gl}(V,\beta)_a = \{x \in \mathrm{End}(V) \mid xV_b \subseteq V_{ab} \text{ for all } b \in G\}$$

and with the bracket

$$[x,y] = xy - \beta(b,a)yx,$$

for all $x \in \mathfrak{gl}(V,\beta)_a, y \in \mathfrak{gl}(V,\beta)_b$. Then $\mathfrak{gl}(V,\beta) = \bigoplus_{a\in G} \mathfrak{gl}(V,\beta)_a$ with this bracket is a Lie color algebra, the so-called *general linear Lie color algebra*.

*The orthosymplectic Lie color algebra $\mathfrak{spo}(V,\beta)$*

A *$\beta$-skew-symmetric bilinear form* is a $\kappa$-bilinear map $\langle\,,\rangle : V \times V \to \kappa$ such that

(1) the form $\langle\,,\rangle$ on $V$ is nondegenerate,

(2) $\langle V_a, V_b\rangle = 0$ if $a \neq b^{-1}$, and

(3) $\langle v,w\rangle = -\beta(b,a)\langle w,v\rangle$, for all $v \in V_a$, $w \in V_b$.

For each $a \in G$, define

$$\mathfrak{spo}(V,\beta)_a = \{x \in \mathfrak{gl}(V,\beta)_a \mid \langle xu,v\rangle + \beta(b,a)\langle u,xv\rangle = 0, \text{ for all } u \in V_b \text{ and } v \in V\}.$$

Then the *orthosymplectic Lie color algebra* is the Lie color subalgebra

$$\mathfrak{spo}(V,\beta) = \bigoplus_{a\in G} \mathfrak{spo}(V,\beta)_a$$

of the Lie color algebra $\mathfrak{gl}(V,\beta)$.

In a similar fashion to Remark 1.1 we define

$$V_{(0)} = \bigoplus_{a,\,\beta(a,a)=1} V_a \quad \text{and} \quad V_{(1)} = \bigoplus_{a,\,\beta(a,a)=-1} V_a .$$

The form $\langle , \rangle$ when restricted to $V_{(0)}$ is skew-symmetric and to $V_{(1)}$ is symmetric. In this way $V$ is a $\mathbb{Z}_2$-graded vector space with a supersymmetric form. Necessarily dim $V_{(0)}$ is even and so throughout this paper we fix

$$\dim\ V_{(0)} = m = 2r, \quad \text{and} \quad n = \dim\ V_{(1)} = n = 2s \quad \text{or} \quad 2s+1.$$

*Remark 1.6.* If $V_{(1)} = (0)$ and $\beta(a,b) = 1$ for all $a,b \in G$, then $\mathfrak{spo}(V,\beta)$ is just the symplectic Lie algebra $\mathfrak{sp}(V_{(0)})$ with a $G$-grading. On the other hand, if $V_{(0)} = (0)$ and $\beta$ is such that

(i) $\beta(b,a) = -1$ for all $a,b$ with $\beta(a,a) = \beta(b,b) = -1$, and
(ii) $\beta(b,a) = 1$ for all $a,b$ such that $\beta(a,a) = 1$ and $\beta(b,b) = -1$,

then $\beta(b,a) = 1$, for all $a,b$ such that $\beta(a,a) = 1 = \beta(b,b)$, and $\mathfrak{spo}(V,\beta)$ is the orthogonal Lie algebra $\mathfrak{so}(V_{(1)})$ with a $G$-grading. In this way the color algebra approach allows us to treat simultaneously the orthogonal and symplectic Lie algebras, the orthosymplectic Lie superalgebra (when $G = \mathbb{Z}_2$), and all the other orthosymplectic color algebras as well.

*Remark 1.7.* It is convenient in what follows to adopt the convention that $\beta(u,x) = \beta(a,b)$ and $\beta(u,v) = \beta(a,b)$ whenever $u \in V_a$, $v \in V_b$, and $x \in gl(V,\beta)_b$ are nonzero. When notation such as $\beta(u,x)$ or $\beta(u,v)$ is used, it is tacitly assumed that the elements are homogeneous.

We assume that the form $\langle , \rangle$ (or the field $\kappa$) is such that there exist homogeneous bases

$$\begin{aligned} B_0 &= \{t_1, t_1^*, t_2, t_2^*, \ldots, t_r, t_r^*\}, \\ B_1 &= \{u_1, u_1^*, u_2, u_2^*, \ldots, u_s, u_s^*, (u_{s+1})\}, \quad \text{and} \\ B &= B_0 \cup B_1, \end{aligned} \tag{1.8}$$

of $V_{(0)}$, $V_{(1)}$, and $V$, respectively, such that

$$\langle v, w\rangle = \langle v^*, w^*\rangle = 0, \quad \text{and} \quad \langle v, w^*\rangle = -\beta(w^*, v)\langle w^*, v\rangle = \delta_{v,w},$$

for all $v, w \in \{t_1, t_2, \ldots, t_r, u_1, u_2, \ldots, u_s, (u_{s+1})\}$. It is to be understood that $u_{s+1}$ occurs only when $n = 2s+1$ and in that case $u_{s+1}^* = u_{s+1}$. We extend the definition of $*$ so that $(t_i^*)^* = t_i$ and $(u_j^*)^* = u_j$ for all $1 \le i \le r$ and all $1 \le j \le s$. Note if $v \in B$ and $v \in V_a$, then $v^* \in V_{a^{-1}}$.

The matrix of the form $\langle\, , \rangle$ relative to the basis $B$ is the matrix $F_B = (F_{b,b'})_{b,b' \in B}$ with $F_{b,b'} = \langle b, b'\rangle$. More explicitly,

$$F_{t_i,t_i^*} = 1, \quad F_{t_i^*,t_i} = -1, \quad F_{u_j,u_j^*} = 1, \quad F_{u_j^*,u_j} = 1, \tag{1.9a}$$

for $i = 1, \ldots, r$, and $j = 1, 2, \ldots, s, (s+1)$, and $F_{b,b'} = 0$ for all other pairs $b, b' \in B$. The inverse of the matrix $F_B$ is the matrix $F_B^{-1} = (F^{-1}_{b,b'})$ with

$$F^{-1}_{t_i,t_i^*} = -1, \quad F^{-1}_{t_i^*,t_i} = 1, \quad F^{-1}_{u_j,u_j^*} = 1, \quad F^{-1}_{u_j^*,u_j} = 1, \tag{1.9b}$$

for $i = 1, \ldots, r$, and $j = 1, 2, \ldots, s, (s+1)$, and $F^{-1}_{b,b'} = 0$ for all other pairs $b, b' \in B$.

*Roots and root vectors in* $\mathfrak{spo}(V, \beta)$

The choice of the basis $B$ of $V$ in (1.8) affords a realization of the elements of $\mathfrak{gl}(V, \beta)$ and $\mathfrak{spo}(V, \beta)$ as matrices. The matrix units $E_{v,w}$, defined by $E_{v,w} y = \delta_{yw} v$ for $v, w, y \in B$, determine a homogeneous basis of $\mathfrak{gl}(V, \beta)$ with $E_{v,w} \in \mathfrak{gl}(V, \beta)_{ab^{-1}}$ whenever $v \in V_a$ and $w \in V_b$.

For $x \in \mathfrak{gl}(V, \beta)$ to belong to $\mathfrak{spo}(V, \beta)$ the relation $\langle xv, w\rangle + \beta(v, x)\langle v, xw\rangle = 0$ must hold, and that translates to the matrix equation

$$v^t x^t F_B w + \beta(v, x) v^t F_B x w = 0,$$

where "t" denotes the usual matrix transpose. When $v, w \in B$, that equation reduces to

$$\begin{aligned} (x^t)_{v,w^*} F_{w^*,w} + \beta(v, x) F_{v,v^*} x_{v^*,w} = 0, \quad \text{or equivalently to} \\ x_{w^*,v^*} F_{w^*,w} + \beta(v^*, x) F_{v^*,v} x_{v,w} = 0, \end{aligned}$$

after replacing $v$ with $v^*$. Using the fact that $F_{v^*,v} = -\beta(v, v^*) F_{v,v*}$ we see that

$$F_{w,w^*} x_{w^*,v^*} = -\beta(v, v)\beta(w, w)\beta(v^*, x) F_{v,v^*} x_{v,w}$$

must hold for all $v, w \in B$. Therefore, whenever $v, w \in \{t_1, \ldots, t_r, u_1, \ldots, u_s, (u_{s+1})\}$, then $x_{w^*,v^*} = -\beta(v, v)\beta(w, w)\beta(v^*, x) x_{v,w}$. For such $v, w$, it follows that $E_{v,w} - \beta(b, b)\beta(a, b) E_{w^*,v^*}$ belongs to $\mathfrak{spo}(V, \beta)_{ab^{-1}}$ if $v \in V_a$ and $w \in V_b$. In particular, the elements

$$h_v = E_{v,v} - E_{v^*,v^*}, \qquad v \in \{t_1, \ldots, t_r, u_1, \ldots, u_s\} \tag{1.10}$$

belong to $\mathfrak{spo}(V, \beta)_{1_G}$, and they span the space $\mathfrak{h}$ of diagonal matrices in $\mathfrak{spo}(V, \beta)$. Now

$$\mathfrak{g} = \bigoplus_{\alpha \in \mathfrak{h}^*} \mathfrak{g}^\alpha, \quad \text{where} \quad \mathfrak{g}^\alpha = \{x \in \mathfrak{g} \mid [h, x] = \alpha(h) x \ \text{ for all } h \in \mathfrak{h}\}.$$

The set $\Delta = \{\alpha \in \mathfrak{h}^* \mid \mathfrak{g}^\alpha \neq (0)\}$ is the set of *roots* of $\mathfrak{g} = \mathfrak{spo}(V, \beta)$ relative to $\mathfrak{h}$. We define

$$\begin{aligned} \Delta_a &= \{\alpha \in \Delta \mid \mathfrak{g}^\alpha \cap \mathfrak{g}_a \neq (0)\}, \quad \text{for each } a \in G, \text{ and} \\ \Delta_0 &= \bigcup_{a,\, \beta(a,a)=1} \Delta_a, \qquad \Delta_1 = \bigcup_{a,\, \beta(a,a)=-1} \Delta_a. \end{aligned}$$

Let $\{\epsilon_v\}$ be the dual basis in $\mathfrak{h}^*$ given by

$$\epsilon_v(h_w) = \delta_{v,w}, \quad \text{for all } v, w \in \{t_1, \cdots, t_r, u_1, \cdots, u_s\},$$

and set $\epsilon_{u_{s+1}} = 0$. Since $E_{v,w}, E_{w^*v^*} \in \mathfrak{gl}(V, \beta)_{ab^{-1}}$ if $v \in V_a$, $w \in V_b$, $v, w \in B$, the elements

$$E_{v,w} - \beta(w, w)\beta(v, w) E_{w^*,v^*} \in \mathfrak{g}^{\epsilon_v - \epsilon_w}, \qquad v, w \in \{t_1, \cdots, t_r, u_1, \cdots, u_s, (u_{s+1})\},\ v \neq w,$$

$$E_{v^*,w} + \beta(v, v)\beta(w, w)\beta(v^*, w) E_{w^*,v} \in \mathfrak{g}^{-\epsilon_v - \epsilon_w}, \qquad v, w \in \{t_1, \cdots, t_r, u_1, \cdots, u_s\}, \text{ where } v \neq w \text{ if } v, w \in \{u_1, \ldots, u_s\}, \tag{1.11}$$

$$E_{v,w^*} + \beta(v^*,w)E_{w,v^*} \in \mathfrak{g}^{\epsilon_v+\epsilon_w}, \qquad \begin{array}{l} v,w \in \{t_1,\cdots,t_r,u_1,\cdots,u_s\}, \\ \text{where } v \neq w \text{ if } v,w \in \{u_1,\ldots,u_s\}, \end{array}$$

along with the elements $h_v$ in (1.10) form a homogeneous basis of $\mathfrak{g} = \mathfrak{spo}(V,\beta)$. Thus, $\dim(\mathfrak{g}^\alpha) = 1$ for all $\alpha \in \Delta$.

As a shorthand set $\epsilon_i = \epsilon_{t_i}$ and $\delta_j = \epsilon_{u_j}$, for $1 \le i \le r$ and $1 \le j \le s$. Then the roots for $\mathfrak{spo}(V,\beta)$ have the following expressions in terms of the $\epsilon_i$'s and $\delta_j$'s:

(i) When $n = 2s+1$,

$\Delta_0 = \{\pm(\epsilon_i \pm \epsilon_j),\ \pm 2\epsilon_i,\ \pm(\delta_k \pm \delta_\ell), \pm\,\delta_k \mid 1 \le i \neq j \le r,\ 1 \le k \neq \ell \le s\}$

$\Delta_1 = \{\pm\epsilon_i,\ \pm(\epsilon_i \pm \delta_j) \mid 1 \le i \le r,\ 1 \le j \le s\}.$

(ii) When $n = 2s$,

$\Delta_0 = \{\pm(\epsilon_i \pm \epsilon_j),\ \pm 2\epsilon_i,\ \pm(\delta_k \pm \delta_\ell) \mid 1 \le i \neq j \le r,\ 1 \le k \neq \ell \le s\}$

$\Delta_1 = \{\pm(\epsilon_i \pm \delta_j) \mid 1 \le i \le r,\ 1 \le j \le s\},$

which are exactly the same as the corresponding sets of even and odd roots for the simple Lie superalgebra $\mathfrak{spo}(V)$ which appears in the classical Lie superalgebra theory see [K2].

If the simple roots are selected to be

$$\alpha_i = \epsilon_i - \epsilon_{i+1}, \quad 1 \le i \le r-1, \qquad \alpha_{r+j} = \delta_j - \delta_{j+1}, \quad 1 \le j \le s-1,$$
$$\alpha_r = \begin{cases} \epsilon_r, & \text{if } n = 1, \\ \epsilon_r - \delta_1, & \text{otherwise}, \end{cases} \qquad \alpha_{r+s} = \begin{cases} \delta_s, & \text{if } n = 2s+1, \\ \delta_{s-1} + \delta_s, & \text{if } n = 2s, \end{cases}$$

then all the simple roots except $\alpha_r$ belong to $\Delta_0$. Every root $\alpha \in \Delta$ is a nonnegative or nonpositive integral combination of the simple roots. Accordingly we write $\alpha > 0$ or $\alpha < 0$. Define elements $x_i \in \mathfrak{g}^{\alpha_i}$, $1 \le i \le r+s$, to be the basis elements in (1.11) given by

$$\begin{aligned} x_i &= E_{t_i,t_{i+1}} - \beta(t_i,t_{i+1})E_{t^*_{i+1},t^*_i}, & & 1 \le i \le r-1, \\ x_r &= E_{t_r,u_1} + \beta(t_r,u_1)E_{u^*_1,t^*_r}, \\ x_{r+j} &= E_{u_j,u_{j+1}} + \beta(u_j,u_{j+1})E_{u^*_{j+1},u^*_j}, & & 1 \le j \le s-1, \\ x_{r+s} &= \begin{cases} E_{u_s,u_{s+1}} + \beta(u_s,u_{s+1})E_{u^*_{s+1},u^*_s}, & \text{if } n = 2s+1, \\ E_{u_s,u^*_{s-1}} + \beta(u^*_s,u_{s-1})E_{u_{s-1},u^*_s}, & \text{if } n = 2s. \end{cases} \end{aligned} \tag{1.12}$$

Then

$$\mathfrak{g} = \mathfrak{g}_- \oplus \mathfrak{h} \oplus \mathfrak{g}_+, \quad \text{where} \quad \mathfrak{g}_- = \bigoplus_{\alpha<0} \mathfrak{g}^\alpha, \quad \text{and} \quad \mathfrak{g}_+ = \bigoplus_{\alpha>0} \mathfrak{g}^\alpha. \tag{1.13}$$

The subalgebra $\mathfrak{g}_+$ in (1.13) is generated by the root vectors $x_i$, $1 \le i \le r+s$.

## 2. THE BRAUER ALGEBRA ACTION ON TENSOR SPACE

*The unfolding map*

A $2k$-*one-factor* is a graph with one row of $2k$ vertices and $k$ edges such that each vertex is incident to precisely one edge. We draw one-factors so that the edges travel above the row of vertices and denote the set of $2k$-one-factors by $\mathcal{O}_{2k}$. We assume that the vertices are numbered 1 to $2k$ from left to right and often represent a one-factor $f \in \mathcal{O}_{2k}$ as a sequence of pairs $f = ((\ell_1, r_1), (\ell_2, r_2), \ldots, (\ell_k, r_k))$, where $\ell_i, r_i \in \{1, \ldots, 2k\}$ give the left vertex and right vertex respectively of each edge. As an example, the sequence of pairs $((1,4)(2,7)(3,5)(6,8))$ represents the 8-one-factor

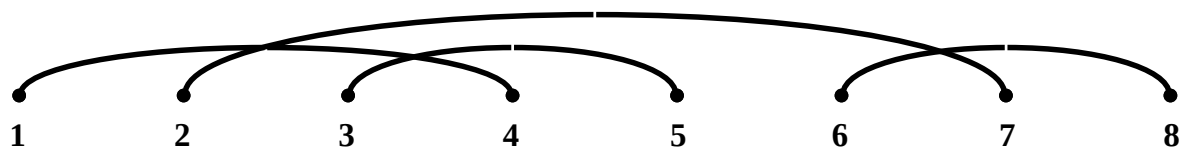

$\in \mathcal{O}_8$.

A *$k$-diagram* is a graph with two rows of $k$ vertices each, one above the other, and $k$ edges such that each vertex is incident to precisely one edge. We draw $k$-diagrams so that the edges remain inside the rectangle formed by the vertices and denote the set of all $k$-diagrams by $\mathcal{D}_k$. We number the top vertices left to right by $1, 2, \ldots, k$ and the bottom vertices left to right by $1', 2', \ldots, k'$.

The *unfolding map*

$$u: \mathcal{D}_k \to \mathcal{O}_{2k} \tag{2.1}$$

converts a $k$-diagram into a $2k$-one-factor by relabeling and repositioning the vertices so that the sequence $(1', 2', \ldots, k', k, (k-1), \ldots, 2, 1)$ becomes $(1, 2, \ldots, 2k)$. As an example, $u((1,4')(2,1')(3,5)(4,6')(6,5')(2',3')) = ((1,11)(2,3)(4,12)(5,7)(6,9)(8,10))$,

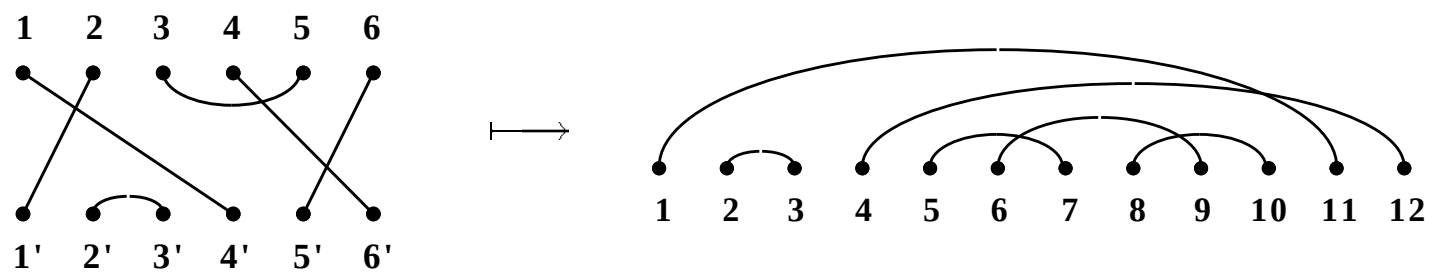

This unfolding map is a combinatorial realization of the isomorphism $\mathrm{End}_{\mathfrak{g}}(V^{\otimes k}) \to (V^{\otimes 2k})^{\mathfrak{g}}$, see Remark 2.13 below.

### *The Brauer algebra*

Let $\kappa$ be a field of characteristic 0 and assume $\eta \in \kappa$. The product of two $k$-diagrams $d_1$ and $d_2$ is obtained by placing $d_1$ above $d_2$ and identifying the vertices in the bottom row of $d_1$ with the corresponding vertices in the top row of $d_2$. The resulting graph contains $k$ paths and some number $c$ of closed loops. If $d$ is the $k$-diagram with the edges that are the paths in this graph but with the closed loops removed, then the product $d_1 d_2$ is given by $d_1 d_2 = \eta^c d$. For example, if

$d_1 =$ 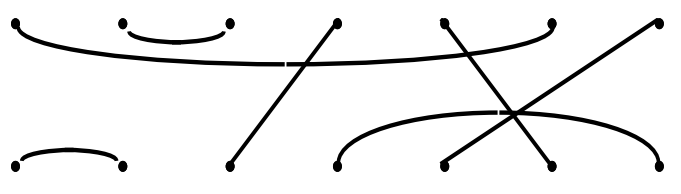 and $d_2 =$ 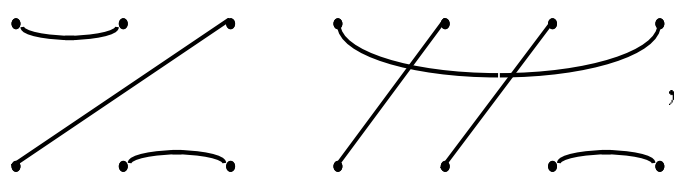,

then

$d_1 d_2 =$ 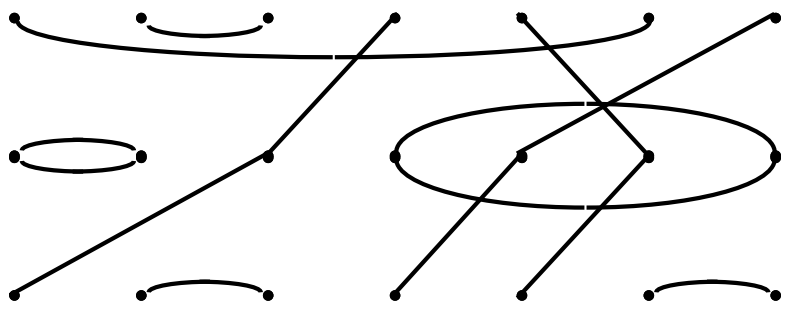 $= \eta^2$ 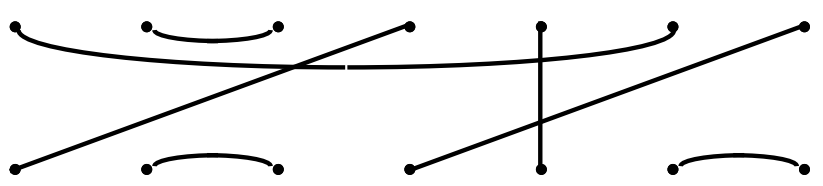.

The *Brauer algebra* $B_k(\eta)$ is the $\kappa$-span of the $k$-diagrams. The $\kappa$-linear extension of the diagram multiplication makes $B_k(\eta)$ into an associative $\kappa$-algebra with identity given by the diagram

$1 =$ 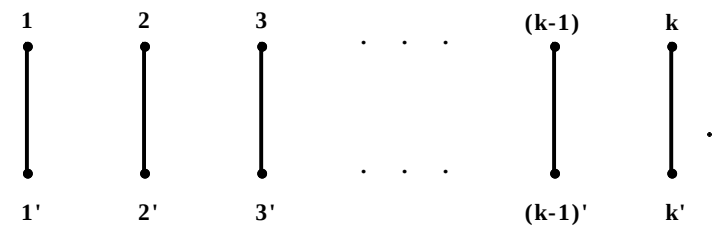

By convention $B_0(\eta) = B_1(\eta) = \kappa$.

The diagrams in $B_k(\eta)$ which have all their edges connecting top vertices to bottom vertices form a symmetric group $S_k$. The elements

$s_i =$ 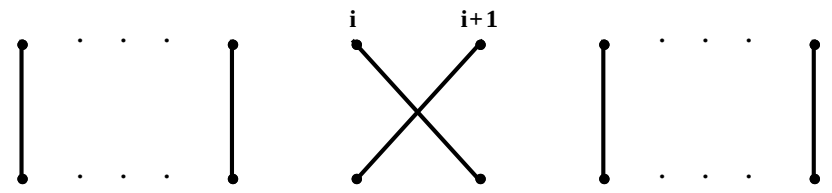

and $e_i =$ 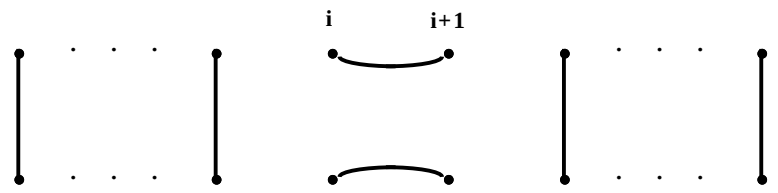

,

$1 \le i \le k-1$, generate the Brauer algebra $B_k(\eta)$. The following theorem is "oft-quoted", (see for example [N]), but we are not aware of a proof in the literature. In the interests of space we shall continue the tradition and not include a proof here either.

**Theorem 2.2.** *The Brauer algebra $B_k(\eta)$ has a presentation as an algebra by generators $s_1, s_2, \ldots, s_{k-1}, e_1, e_2, \ldots, e_{k-1}$ and relations*

$$s_i^2 = 1, \quad e_i^2 = \eta e_i, \quad e_i s_i = s_i e_i = e_i, \qquad 1 \le i \le k-1,$$

$$s_i s_j = s_j s_i, \quad s_i e_j = e_j s_i, \quad e_i e_j = e_j e_i, \qquad |i-j| > 1,$$

$$s_i s_{i+1} s_i = s_{i+1} s_i s_{i+1}, \quad e_i e_{i+1} e_i = e_i, \quad e_{i+1} e_i e_{i+1} = e_{i+1}, \qquad 1 \le i \le k-2,$$

$$s_i e_{i+1} e_i = s_{i+1} e_i, \quad e_{i+1} e_i s_{i+1} = e_{i+1} s_i, \qquad 1 \le i \le k-2.$$

*Some facts about braided rigid monoidal categories*

Let $\mathcal{C}$ be a $\kappa$-linear braided rigid monoidal category whose identity object we denote by $\mathbf{1}$. The defining property of $\mathbf{1}$ is that there are natural isomorphisms $V \otimes \mathbf{1} \cong V \cong \mathbf{1} \otimes V$ for all $V \in \mathcal{C}$. For simplicity we assume that $\mathcal{C}$ is strict. Although everything we do in this section works in complete generality for any $\kappa$-linear braided rigid monoidal category that is strict, the reader who prefers not to work with general categorical language can just note that for the category of finite dimensional modules for a Lie color algebra $\mathfrak{g}$, we have the following:

(2.3) Let $V$ be an object in $\mathcal{C}$. The space of *invariants* in $V$ is the $\kappa$-vector space defined by

$$V^{\mathcal{C}} = \mathrm{Hom}_{\mathcal{C}}(\mathbf{1}, V).$$

(2.4) Every object $V \in \mathcal{C}$ has a *dual*, i.e. there is an object $V^*$ in $\mathcal{C}$ and morphisms

$$\mathrm{pr}_V : \mathbf{1} \to V \otimes V^*, \qquad \text{and} \qquad \mathrm{ev}_V : V^* \otimes V \to \mathbf{1}.$$

The dual of $V$ is unique up to isomorphism.

(2.5) If $M$ and $N$ are any two objects in $\mathcal{C}$, then

$$(M \otimes N)^* \cong N^* \otimes M^*;$$

in fact, this isomorphism can be explicitly seen using the canonical maps

$$\mathbf{1} \xrightarrow{\mathrm{pr}_M} M \otimes M^* \cong M \otimes \mathbf{1} \otimes M^* \xrightarrow{\mathrm{id} \otimes \mathrm{pr}_N \otimes \mathrm{id}} M \otimes N \otimes N^* \otimes M^*$$

and

$$N^* \otimes M^* \otimes M \otimes N \xrightarrow{\mathrm{id} \otimes \mathrm{ev}_M \otimes \mathrm{id}} N^* \otimes \mathbf{1} \otimes N \cong N^* \otimes N \xrightarrow{\mathrm{ev}_N} \mathbf{1}.$$

(2.6) There is an isomorphism of $\kappa$-vector spaces

$$\begin{array}{ccc} \mathrm{Hom}_{\mathcal{C}}(M, P \otimes N^*) & \to & \mathrm{Hom}_{\mathcal{C}}(M \otimes N, P) \\ \psi & \mapsto & \phi \end{array},$$

where $\phi$ is the morphism defined by

$$M \otimes N \xrightarrow{\psi \otimes \mathrm{id}} P \otimes N^* \otimes N \xrightarrow{\mathrm{id} \otimes \mathrm{ev}_N} P \otimes \mathbf{1} \cong P.$$

The inverse map is

$$\begin{array}{ccc} \mathrm{Hom}_{\mathcal{C}}(M \otimes N, P) & \to & \mathrm{Hom}_{\mathcal{C}}(M, P \otimes N^*) \\ \phi & \mapsto & \psi \end{array},$$

where $\psi$ is given by

$$M \cong M \otimes \mathbf{1} \xrightarrow{\mathrm{id} \otimes \mathrm{pr}_N} M \otimes N \otimes N^* \xrightarrow{\phi \otimes \mathrm{id}} P \otimes N^*.$$

(2.7) If $U$ and $V$ are any two objects in $\mathcal{C}$, then there is an isomorphism

$$\check{R}_{U,V} : U \otimes V \to V \otimes U.$$

(2.8) If $V$ is an object in $\mathcal{C}$, then the isomorphisms

$$\check{R}_i = \mathrm{id}^{\otimes(i-1)} \otimes (-\check{R}_{V,V}) \otimes \mathrm{id}^{\otimes(k-i-1)}$$

satisfy the properties

$$\begin{aligned} \check{R}_i \check{R}_j &= \check{R}_j \check{R}_i, \qquad \text{if } |i-j| \geq 2, \\ \check{R}_i \check{R}_{i+1} \check{R}_i &= \check{R}_{i+1} \check{R}_i \check{R}_{i+1}. \end{aligned} \tag{2.9}$$

*The braid group action on* $V^{\otimes k}$

The relations in (2.9) imply that the $R_i$'s afford a representation of the braid group on $V^{\otimes k}$, (see for example, [CP, §15.2A]). Note that since $\check{R}_{V,V}$ is an isomorphism in $\mathrm{Hom}_{\mathcal{C}}(V^{\otimes 2}, V^{\otimes 2})$, so is $-\check{R}_{V,V}$. The only reason that we are adjusting by a negative sign is so that the $\check{R}_i$ will be the images of the Brauer algebra generators in the map $\Psi$ which will be discussed in Theorem 2.16.

If $\pi$ is a permutation in the symmetric group define

$$\check{R}_\pi = \check{R}_{i_1}\check{R}_{i_2}\cdots\check{R}_{i_p},$$

where $\pi = s_{i_1}s_{i_2}\cdots s_{i_p}$ is a reduced expression of $\pi$ (a reduced expression for $\pi$ is a product of simple reflections $s_i = (i,i+1)$ with $p$ as small as possible). Since the relations in (2.9) hold, the map $\check{R}_\pi \in \mathrm{Hom}_{\mathcal{C}}(V^{\otimes k}, V^{\otimes k})$ is independent of the choice of reduced expression (see for example, [Bo, Prop. 5, Chap. IV, §1, no.5]). If $f \in \mathcal{O}_{2k}$ is a $2k$-one-factor, say $f = ((\ell_1,r_1)(\ell_2,r_2)\cdots(\ell_k,r_k))$ with $\ell_1 < \ell_2 < \cdots < \ell_k$, let $\pi_f$ be the permutation in the symmetric group $S_{2k}$ given by

$$\pi_f = \begin{pmatrix} 1 & 2 & \cdots & k & k+1 & k+2 & \cdots & 2k \\ \ell_1 & \ell_2 & \cdots & \ell_k & r_k & r_{k-1} & \cdots & r_1 \end{pmatrix} \in S_{2k}.$$

Then define

$$\check{R}_f = \check{R}_{\pi_f}, \quad \text{for } f \in \mathcal{O}_{2k}, \qquad \text{and} \qquad \check{R}_d = \check{R}_{u(d)}, \quad \text{for } d \in \mathcal{D}_k, \tag{2.10}$$

where $u\colon \mathcal{D}_k \to \mathcal{O}_{2k}$ is the unfolding map of (2.1).

*The maps $\Phi$ and $\Psi$*

Let $V$ be an object of $\mathcal{C}$ such that $V \cong V^*$. Fix an isomorphism

$$F\colon V \xrightarrow{\sim} V^*$$

in $\mathrm{Hom}_{\mathcal{C}}(V,V^*)$. From (2.5) there is a natural isomorphism $(V^{\otimes k})^* \cong (V^*)^{\otimes k}$. Let $\widetilde{\mathrm{pr}}_k$ and $\widetilde{\mathrm{ev}}_k$ be the compositions

$$\begin{aligned} \widetilde{\mathrm{pr}}_k &: \mathbf{1} \xrightarrow{\mathrm{pr}_{V^{\otimes k}}} V^{\otimes k} \otimes (V^{\otimes k})^* \xrightarrow{\sim} V^{\otimes k} \otimes (V^*)^{\otimes k}, \qquad \text{and} \\ \widetilde{\mathrm{ev}}_k &: (V^*)^{\otimes k} \otimes V^{\otimes k} \xrightarrow{\sim} (V^{\otimes k})^* \otimes V^{\otimes k} \xrightarrow{\mathrm{ev}_{V^{\otimes k}}} \mathbf{1}. \end{aligned}$$

Now define a map

$$\begin{array}{cccc} \Phi\colon & \mathcal{O}_{2k} & \longrightarrow & \mathrm{Hom}_{\mathcal{C}}(\mathbf{1}, V^{\otimes 2k}) = (V^{\otimes 2k})^{\mathcal{C}} \\ & f & \longmapsto & \Phi_f \end{array}$$

by letting $\Phi_f$ be the composition

$$\Phi_f : \mathbf{1} \xrightarrow{\widetilde{\mathrm{pr}}_k} V^{\otimes k} \otimes (V^*)^{\otimes k} \xrightarrow{\mathrm{id}^{\otimes k} \otimes (F^{-1})^{\otimes k}} V^{\otimes 2k} \xrightarrow{\check{R}_f} V^{\otimes 2k}. \tag{2.11}$$

Similarly, define a map

$$\begin{array}{cccc} \Psi\colon & \mathcal{D}_k & \longrightarrow & \mathrm{Hom}_{\mathcal{C}}(V^{\otimes k}, V^{\otimes k}) \\ & d & \longmapsto & \Psi_d \end{array}$$

by setting $\Psi_d$ equal to the composite map

$$\begin{aligned} \Psi_d : V^{\otimes k} \cong \mathbf{1} \otimes V^{\otimes k} &\xrightarrow{\widetilde{\mathrm{pr}}_k \otimes \mathrm{id}^{\otimes k}} V^{\otimes k} \otimes (V^*)^{\otimes k} \otimes V^{\otimes k} \xrightarrow{\mathrm{id}^{\otimes k} \otimes (F^{-1})^{\otimes k} \otimes \mathrm{id}^{\otimes k}} V^{\otimes 2k} \otimes V^{\otimes k} \\ &\xrightarrow{\check{R}_d \otimes \mathrm{id}^{\otimes k}} V^{\otimes 2k} \otimes V^{\otimes k} \xrightarrow{\mathrm{id}^{\otimes k} \otimes F^{\otimes k} \otimes \mathrm{id}^{\otimes k}} V^{\otimes k} \otimes (V^*)^{\otimes k} \otimes V^{\otimes k} \\ &\xrightarrow{\mathrm{id}^{\otimes k} \otimes \widetilde{\mathrm{ev}}_k} V^{\otimes k} \otimes \mathbf{1} \cong V^{\otimes k}. \end{aligned} \tag{2.12}$$

The map $\Phi$ associates an invariant $\Phi_f$ in $(V^{\otimes k})^{\mathcal{C}}$ to each $2k$-one-factor $f$. The map $\Psi$ allows us to define an action of each $k$-diagram $d$ on the vector space $V^{\otimes k}$ by letting $d$ act by the homomorphism $\Psi_d$.

*Remark 2.13.* The maps $\Phi$ and $\Psi$ are related. Indeed, $\Psi_d$ is the image of $(\mathrm{id}^{\otimes k} \otimes F^{\otimes k}) \circ \Phi_{u(d)} \in \mathrm{Hom}_{\mathcal{C}}(\mathbf{1}, V^{\otimes k} \otimes (V^{\otimes k})^*)$ under the isomorphism $\mathrm{Hom}_{\mathcal{C}}(\mathbf{1}, V^{\otimes k} \otimes (V^{\otimes k})^*) \cong \mathrm{Hom}_{\mathcal{C}}(V^{\otimes k}, V^{\otimes k})$ which is given by setting $M = \mathbf{1}$ and $N = P = V^{\otimes k}$ in (2.6).

*The commuting action of the Brauer algebra on the $\mathfrak{spo}(V,\beta)$-module $V^{\otimes k}$*

Let $G$ be a finite abelian group, $\beta$ a symmetric bicharacter on $G$, $V$ a $G$-graded vector space with a nondegenerate $\beta$-skew-symmetric bilinear form $\langle , \rangle$ and $\mathfrak{g} = \mathfrak{spo}(V,\beta)$. There is a $\mathfrak{g}$-module isomorphism

$$\begin{array}{rccc} F\colon & V & \longrightarrow & V^* \\ & v & \longmapsto & \langle v, \cdot \rangle. \end{array}$$

In this special case, we will compute the map $\Psi$ from (2.12) explicitly.

Let $B = \{v_1, \ldots, v_N\}$ be a homogeneous basis of $V$, i.e. for each $1 \le i \le N$, $v_i \in V_{c_i}$ for some $c_i \in G$. Let $\{v^1, \ldots, v^N\}$ be the dual basis in $V^*$. Assume

$$F_B = (F_{i,j})_{1 \le i,j \le N}, \qquad \text{where} \qquad \langle v_i, v_j \rangle = F_{i,j}, \tag{2.14}$$

is the matrix of the form $\langle , \rangle$ with respect to the basis $B$. Sometimes we will write $F_{v_i,v_j}$ instead of $F_{i,j}$ as we have done in the previous section. Let $F_B^{-1} = (F_{i,j}^{-1})_{1 \le i,j \le N}$ be the inverse of the matrix $F_B$. Then

$$F(v_i) = \sum_{j=1}^{N} F_{i,j} v^j, \qquad \text{and} \qquad F^{-1}(v^j) = \sum_{i=1}^{N} F_{j,i}^{-1} v_i. \tag{2.15}$$

Let $d$ be a $k$-diagram. Label the top vertices (left to right) with a sequence $\underline{a} = (a_1, a_2, \ldots, a_k)$ of basis elements $a_i \in B$ and the bottom vertices (left to right) with a sequence $\underline{b} = (b_1, b_2, \ldots, b_k)$ of basis elements $b_i \in B$. Assign a weight to each edge and each crossing of this labeled $k$-diagram according to the following:

(1) If a horizontal edge $(a, a')$ on the top has $a$ to the left of $a'$, assign the weight $F_{a,a'}$ to it.
(2) If a horizontal edge $(b, b')$ on the bottom has $b$ to the left of $b'$, weight it by $F_{b,b'}^{-1}$.
(3) Weight each vertical edge $(a, b)$ by $\delta_{a,b}$ (Kronecker delta),
(4) Weight each crossing by $-\beta(\ell_1, \ell_2)$, where $\ell_1$ is the label of a vertex adjacent to the first edge, and $\ell_2$ is a vertex adjacent to the second edge in the crossing. Of the four vertices adjacent to the two edges that cross, $\ell_1$ and $\ell_2$ should be chosen to be the last two vertices (in order) when counting off the vertices in a counterclockwise fashion beginning from the bottom left corner of the diagram.

The *weight* of the labeled $k$-diagram, which we denote $d_{\underline{a},\underline{b}}$, is the product of the weights over all the edges and crossings. For example, the weight of the following labeled 7-diagram is the product $(-1)^5 F_{3,3}F_{1,2}F_{2,2}^{-1}F_{2,4}^{-1}\delta_{1,2}\delta_{1,3}\delta_{2,4}\beta(v_2,v_3)\beta(v_1,v_3)\beta(v_4,v_1)^2\beta(v_4,v_4) = 0$. We have adopted an abbreviated notation in displaying the weights of the edges so as not to clutter the picture.

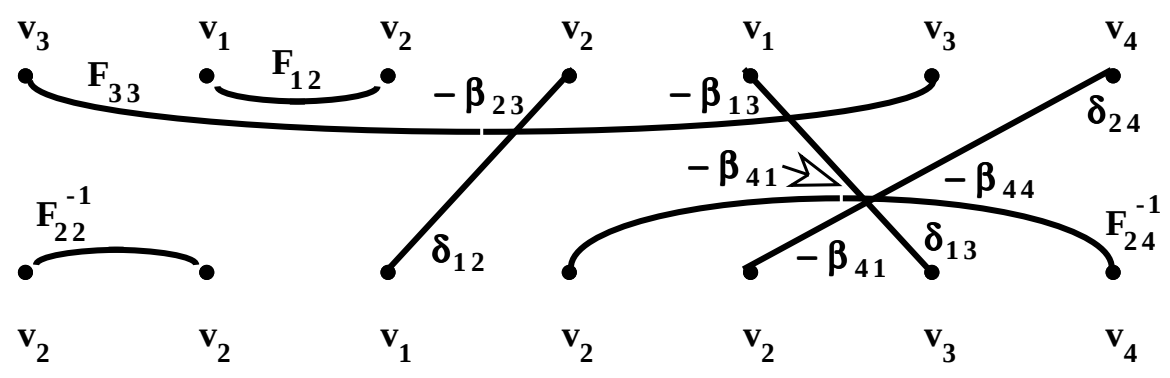

**Theorem 2.16.** *Let $G$ be a finite group and let $\beta$ be a symmetric bicharacter on $G$. Assume $V$ is a $G$-graded vector space with a $\beta$-skew-symmetric bilinear form $\langle , \rangle$. Let $\mathfrak{g} = \mathfrak{spo}(V,\beta)$ be the orthosymplectic Lie color algebra.*

(a) *The image of a $k$-diagram $d \in \mathcal{D}_k$, under the map $\Psi\colon \mathcal{D}_k \to \mathrm{Hom}_{\mathfrak{g}}(V^{\otimes k}, V^{\otimes k})$ from (2.12), is the homomorphism $\Psi_d$ given explicitly by*

$$(a_1 \otimes \cdots \otimes a_k)\Psi_d = \sum_{b_1,\ldots,b_k \in B} d_{\underline{a},\underline{b}} b_1 \otimes \cdots \otimes b_k,$$

*where $d_{\underline{a},\underline{b}}$ is the weight of the $k$-diagram $d$ with top vertices labeled by $a_1,\ldots,a_k$ and bottom vertices labeled by $b_1,\ldots,b_k$. (We have chosen to write $\Psi_d$ as an operator on the right so that (b) will hold.)*

(b) *Let $m = \dim V_{(0)}$ and $n = \dim V_{(1)}$ where $V_{(0)} = \sum_{c\in G,\beta(c,c)=1} V_c$ and $V_{(1)} = \sum_{c\in G,\beta(c,c)=-1} V_c$. Then the map $\Psi$ extends to a homomorphism*

$$\Psi\colon B_k(n-m) \to \mathrm{Hom}_{\mathfrak{g}}(V^{\otimes k}, V^{\otimes k})$$

*of algebras.*

*Proof.* (a) In this setting the map $\widetilde{\mathrm{pr}}_k$ is given by

$$\begin{aligned} \widetilde{\mathrm{pr}}_k\colon \mathbf{1} &\longrightarrow V^{\otimes k} \otimes (V^*)^{\otimes k} \\ 1 &\longmapsto \sum_{1\le i_1,\ldots,i_k\le N} v_{i_1}\otimes\cdots\otimes v_{i_k}\otimes v^{i_k}\otimes\cdots\otimes v^{i_1}, \end{aligned}$$

which can be seen by induction using (2.5). One computes $\Psi_d$ from the definition in (2.12). The following sequence gives a representative pictorial example, but the general case is done in exactly the same fashion. For brevity in the pictures we reduce the notation even further and write $i_3$ for $v_{i_3}$ and $i^3$ for the dual vector $v^{i_3}$.

Let $a_i = v_{n_i}$ for $i = 1,2,3,4$ and let us compute the action of $\Psi_d$ for the diagram

$d =$ 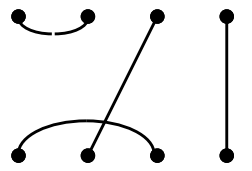

on the element $a_1a_2a_3a_4 = v_{n_1}v_{n_2}v_{n_3}v_{n_4} \in V^{\otimes 4}$ (we omit the tensor symbols $\otimes$ as another space-saving device). Note that the image of the diagram $d$ under the unfolding map $u$ of (2.1) is the one-factor

$u(d) =$ 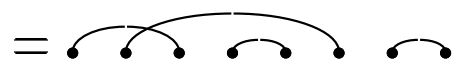 .

Following the definition of $\Psi_d$ as a composite map and using the fact that $\check{R}_{V,V}(v \otimes w) =$

$-\beta(w,v)w\otimes v$, we have

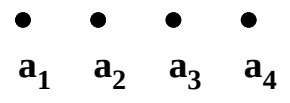

$\xrightarrow{\tilde{pr}_k\otimes \mathrm{id}^{\otimes k}} \sum_{i_1,i_2,i_3,i_4}$

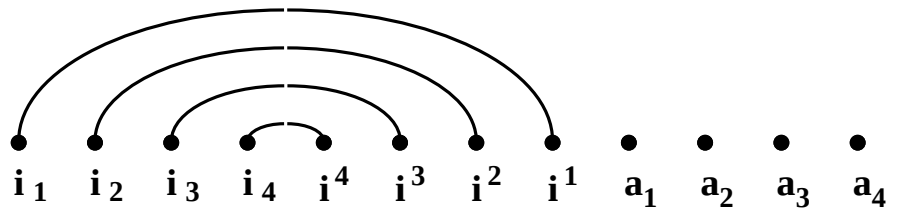

$\xrightarrow{\mathrm{id}^{\otimes k}\otimes(F^{-1})^{\otimes k}\otimes \mathrm{id}^{\otimes k}} \sum_{i_1,i_2,i_3,i_4}$

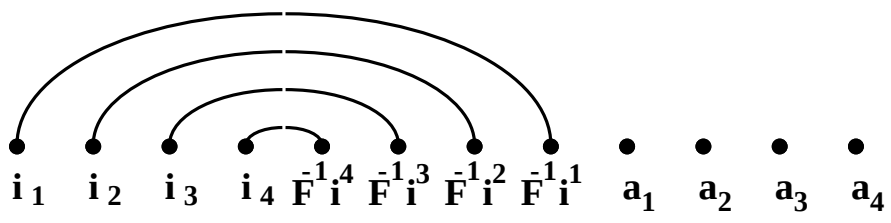

$\xrightarrow{\check{R}_d\otimes \mathrm{id}^{\otimes k}} \sum_{i_1,i_2,i_3,i_4} -\beta(F^{-1}v^{i_1},F^{-1}v^{i_2})$

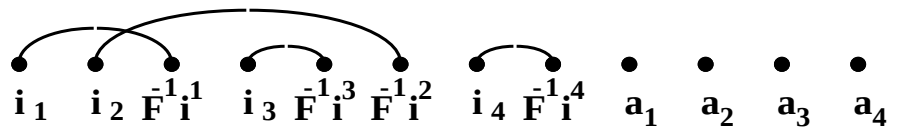

$\xrightarrow{\mathrm{id}^{\otimes k}\otimes F^{\otimes k}\otimes \mathrm{id}^{\otimes k}} \sum_{i_1,i_2,i_3,i_4} -\beta(F^{-1}v^{i_1},v^{i_2})$

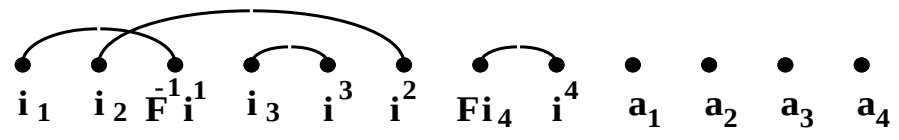

$\xrightarrow{\mathrm{id}^{\otimes k}\otimes \tilde{ev}_k} \sum_{i_1,i_2,i_3,i_4} -\beta(F^{-1}v^{i_1},v^{i_2})$

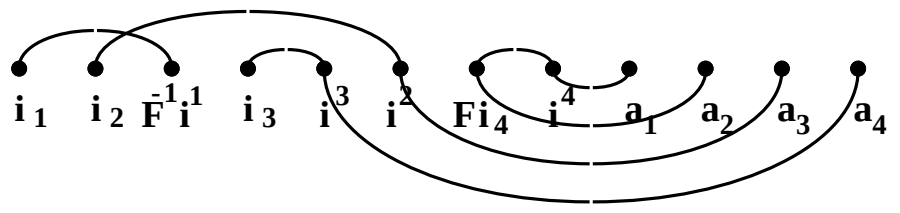

$$= \sum_{i_1,i_2,i_3,i_4,j_1,j_4} -\beta(v_{j_1},v^{i_2})F^{-1}_{i_1,j_1}\langle v^{i_3},v_{n_4}\rangle\langle v^{i_2},v_{n_3}\rangle\langle v^{j_4},v_{n_2}\rangle\langle v^{i_4},v_{n_1}\rangle F_{i_4,j_4}v_{i_1}v_{i_2}v_{j_1}v_{i_3}$$

$$= \sum_{i_1,i_2,j_1,i_3} \mathrm{wt}\left(\quad\right) v_{i_1}v_{i_2}v_{j_1}v_{i_3}.$$

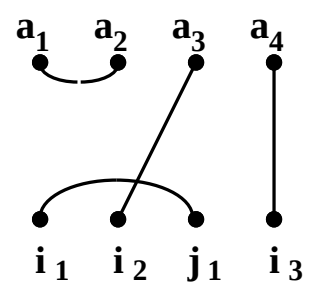

In doing these calculations it is helpful to note that $\deg(Fv_i)=\deg(v_i)=\deg(v^i)^{-1}=\deg(F^{-1}v^i)^{-1}$.

(b) It follows from the definition of the weights of labeled diagrams that

$$(v_{i_1}v_{i_2})\Psi_e = F_{i_1,i_2}\sum_{j_1,j_2}F^{-1}_{j_1,j_2}v_{j_1}v_{j_2},\quad \text{and}\quad (v_{i_1}v_{i_2})\Psi_s = -\beta(v_{i_2},v_{i_1})v_{i_2}v_{i_1}, \tag{2.17}$$

for the case $k=2$ and the 2-diagrams

$e =$ 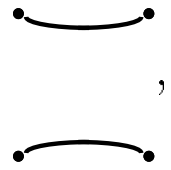 , and $s =$ 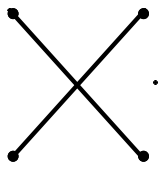 .

For arbitrary $k$,

$$\Psi_{s_i} = \mathrm{id}^{\otimes(i-1)} \otimes \Psi_s \otimes \mathrm{id}^{\otimes(k-i-1)}, \qquad \text{and}$$
$$\Psi_{e_i} = \mathrm{id}^{\otimes(i-1)} \otimes \Psi_e \otimes \mathrm{id}^{\otimes(k-i-1)},$$

where $e_i$ and $s_i$ are the generators of the Brauer algebra in Theorem 2.2. Given these observations it is sufficient to check that $\Psi_{e_i}$ and $\Psi_{s_i}$ satisfy the relations in the statement of Theorem 2.2, for $i = 1$. All of these relations are easily verified by direct calculation. For example, the relation $e_1^2 = (n-m)e_1$ follows from the computation

$$\begin{aligned}
(v_{i_1}v_{i_2})\Psi_{e_1}^2 &= F_{i_1,i_2}\sum_{\ell_1,\ell_2} v_{\ell_1}v_{\ell_2}F^{-1}_{\ell_1,\ell_2}\sum_{j_1,j_2} F_{j_1,j_2}F^{-1}_{j_1,j_2} \\
&= F_{i_1,i_2}\sum_{\ell_1,\ell_2} v_{\ell_1}v_{\ell_2}F^{-1}_{\ell_1,\ell_2}\sum_{j_1,j_2} -\beta(v_{j_2},v_{j_1})F_{j_2,j_1}F^{-1}_{j_1,j_2} \\
&= F_{i_1,i_2}\sum_{\ell_1,\ell_2} v_{\ell_1}v_{\ell_2}F^{-1}_{\ell_1,\ell_2}\sum_{j_2} -\beta(v_{j_2},v_{j_2})\sum_{j_1} F_{j_2,j_1}F^{-1}_{j_1,j_2} \\
&= F_{i_1,i_2}\sum_{\ell_1,\ell_2} v_{\ell_1}v_{\ell_2}F^{-1}_{\ell_1,\ell_2}\sum_{j_2} -\beta(v_{j_2},v_{j_2})\delta_{j_2,j_2} \\
&= (n-m)F_{i_1,i_2}\sum_{\ell_1,\ell_2} v_{\ell_1}v_{\ell_2}F^{-1}_{\ell_1,\ell_2} \\
&= (v_{i_1}v_{i_2})(n-m)\Psi_{e_1}.
\end{aligned}$$

In this argument the critical facts used are that $F_{j_1,j_2} = -\beta(v_{j_2},v_{j_1})F_{j_2,j_1}$ and that $F_{j_2,j_1} \neq 0$ forces $\deg(v_{j_1}) = \deg(v_{j_2})^{-1}$ so that $\beta(v_{j_2},v_{j_1}) = \beta(v_{j_2},v_{j_2})^{-1} = \beta(v_{j_2},v_{j_2}) = \pm 1$.

The relation $\Psi_{s_1}^2 = \mathrm{id}^{\otimes k}$ follows from the fact that $\beta(i_1,i_2)\beta(i_2,i_1) = 1$, since $\beta$ is symmetric. The relation $\Psi_{e_1}\Psi_{s_1} = \Psi_{s_1}\Psi_{e_1} = \Psi_{e_1}$ can be derived from the same type of easy calculations. One shows that

$$(v_{i_1}v_{i_2}v_{i_3})\Psi_{s_1}\Psi_{s_2}\Psi_{s_1} = (v_{i_1}v_{i_2}v_{i_3})\Psi_{s_2}\Psi_{s_1}\Psi_{s_2} = (-\beta(i_3,i_1))(-\beta(i_2,i_1))(-\beta(i_3,i_2))v_{i_3}v_{i_2}v_{i_1}$$

and the remainder of the relations are checked with similar calculations. ∎

### $\mathfrak{spo}(V,\beta)$*-invariants in* $V^{\otimes 2k}$

We retain the setup and notation of Theorem 2.16. If $f$ is a $2k$-one-factor, let $f_{\underline{b}}$ denote the $2k$-one-factor $f$ with its vertices labeled (left to right) by the sequence $\underline{b} = (b_1,\ldots,b_{2k})$ of basis elements $b_1,\ldots,b_{2k} \in B$. Assign a weight to each edge and each crossing of the labeled one-factor $f_{\underline{b}}$ as follows:

(a) Weight the edge $(b,b')$ by $F^{-1}_{b,b'}$.

(b) Weight each crossing by $-\beta(\ell_1,\ell_2)$ where $\ell_1$ is the label of a vertex that is adjacent to one of the edges, and $\ell_2$ is the label of a vertex that is adjacent to the other vertex. Of the four vertices adjacent to the two edges that cross; $\ell_1$ and $\ell_2$ should be chosen to be the last two vertices (in order) when counting off the vertices left to right.

Given a one-factor $f_{\underline{b}}$ which is labeled by $\underline{b} = (b_1,\ldots,b_{2k})$ define the weight $\mathrm{wt}(f_{\underline{b}})$ of $f_{\underline{b}}$ to be the product over all the weights of the crossings and edges. The proof of the following proposition is exactly analogous to the proof of part (a) of Theorem 2.16.

**Proposition 2.18.** *Retain the setup and notation of Theorem 2.16. The image of a $2k$-one-factor $f \in \mathcal{O}_{2k}$, under the map $\Phi : \mathcal{O}_{2k} \to (V^{\otimes 2k})^{\mathfrak{g}}$ from (2.11), is the invariant $\Phi_f$ given explicitly by*

$$\Phi_f = \Phi_f(1) = \sum_{b_1,\dots,b_{2k}\in B} \mathrm{wt}(f_{\underline{b}}) b_1 \otimes \cdots \otimes b_{2k},$$

*where the sum is over all possible labelings $\underline{b} = (b_1, b_2, \dots, b_{2k})$ of the one-factor $f$.*

## 3. MAXIMAL VECTORS IN TENSOR SPACE

Let $\mathfrak{h}$ denote the diagonal matrices in $\mathfrak{spo}(V,\beta)$ as in Section 1. An $\mathfrak{spo}(V,\beta)$ -module $M$ is said to have a weight space decomposition relative to $\mathfrak{h}$ if $M = \oplus \sum_{\mu\in\mathfrak{h}^*} M_\mu$, where $M_\mu = \{v \in M \mid hv = \mu(h)v \text{ for all } h \in \mathfrak{h}\}$. A vector $v_+$ in $M$ is a maximal vector of weight $\lambda$ if $\mathfrak{g}_+ v_+ = (0)$ and $hv_+ = \lambda(h)v_+$ for all $h \in \mathfrak{h}$. In particular, the module $V$ is an irreducible $\mathfrak{spo}(V,\beta)$- module with a unique (up to scalar multiple) maximal vector $v_+ = t_1$, which has weight $\epsilon_1$. All the basis vectors in $B = \{t_1, t_1^*, \dots, t_r, t_r^*, u_1, u_1^*, \dots, u_s, u_s^*, (u_{s+1})\}$ are weight vectors with weights given by

$$\begin{aligned}
&\mathrm{wt}(t_i) = \epsilon_i = -\mathrm{wt}(t_i^*), && i = 1,\dots,r,\\
&\mathrm{wt}(u_j) = \delta_j = -\mathrm{wt}(u_j^*), && j = 1,\dots,s,\\
&\mathrm{wt}(u_{s+1}) = 0, && \text{if } n = 2s+1.
\end{aligned} \tag{3.1}$$

If $v_1 \otimes \cdots \otimes v_k$ is a simple tensor in $V^{\otimes k}$ with $v_i \in B$ for each $i$, then its weight is just the sum of the weights of its factors.

In this section we construct maximal vectors of $V^{\otimes k}$ using the elements $s_i$, $e_i$ coming from the Brauer algebra. The maximal vectors in $V^{\otimes k}$ are often a help in locating its irreducible summands. We begin with a few combinatorial preliminaries.

### *Young symmetrizers, contractions, and $\mathfrak{spo}(V,\beta)$-submodules of $V^{\otimes k}$*

Suppose $\lambda$ is a partition of $k$. The *standard tableaux* $T$ of shape $\lambda$ are obtained by filling in the frame of $\lambda$ with the elements of $\mathcal{K} = \{1,\dots,k\}$ so that the entries increase across the rows from left to right and down the columns. Let $\mathcal{ST}(\mathcal{K})$ denote the set of all standard tableaux as $\lambda$ ranges over all the partitions of $k$. We associate two subgroups in the symmetric group $S_k$ to each standard tableau $T \in \mathcal{ST}(\mathcal{K})$. The first is the *row group* $R_T$ consisting of all permutations in $S_k$ which permute only the elements of $\mathcal{K}$ within each row of $T$, while the second is the *column group* $C_T$ of $T$ or group of permutations moving only the elements of $\mathcal{K}$ in each column. Then the element

$$s_T = \left(\sum_{\psi\in R_T} sgn(\psi)\psi\right)\left(\sum_{\gamma\in C_T}\gamma\right) = \sum_{\substack{\psi\in R_T\\ \gamma\in C_T}} sgn(\psi)\psi\gamma \in \mathbb{C}S_k \tag{3.2}$$

has the property that there is some $h(\lambda) \in \mathbb{Z}^+$ that only depends on the underlying frame of $T$ such that $s_T^2 = h(\lambda)s_T$ (see [We, Chap. 4, §2]). If $y_T = (1/h(\lambda))s_T$, then $y_T^2 = y_T$ so that $y_T$ is idempotent. We refer to $y_T$ as the *symmetrizer* determined by $T$. For each $T \in \mathcal{ST}(\mathcal{K})$ the space $\mathbb{C}S_k y_T$ is a minimal left ideal of the group algebra $\mathbb{C}S_k$ and $\mathbb{C}S_k = \oplus \sum_{T\in\mathcal{ST}(\mathcal{K})} \mathbb{C}S_k y_T$, (see for example, [We, Chap. 4, §3.4]).

Starting at the left end of the first row and moving from left to right we compare the entries of two standard tableaux $T_1$ and $T_2$ of shape $\lambda$. If the first nonzero difference $j_1 - j_2$ is positive for corresponding entries $j_1$ in $T_1$ and $j_2$ in $T_2$, then we say $T_1 > T_2$. If all corresponding entries in the first row are equal, then we proceed to the second row, etc. Thus, $T_1 > T_2$ if $j_1 > j_2$ for the first pair of corresponding entries $j_1$, $j_2$ which differ. With respect to this ordering, we have the following lemma.

**Lemma 3.3.** *(See for example, [We, Theorem 4.3D].) Let $T_1$, $T_2 \in \mathcal{ST}(\mathcal{K})$ be two standard tableaux of the same shape. If $T_1 > T_2$, then $y_{T_1} y_{T_2} = 0$.*

A partition $\lambda = (\lambda_1, \ldots, \lambda_\ell, 0, \ldots) \vdash l$ is said to be of $(r,s)$*-hook shape* if $\lambda_{r+1} \leq s$. If $\mathcal{L}$ is a subset of $\mathcal{K}$ of cardinality $l$, then

$$\mathcal{HST}_{r,s}(\mathcal{L}) \stackrel{\text{def}}{=} \{T \in \mathcal{ST}(\mathcal{L}) \mid T \text{ has an } (r,s)\text{–hook shape}\} \tag{3.4}$$

is the set of all standard tableaux with entries in $\mathcal{L}$ of shape $\lambda$ as $\lambda$ ranges over all partitions of $l$ that are $(r,s)$-hook shape.

**Example.** Suppose $r = 2$ and $s \geq 2$, and let $\lambda = (1^2, 2^2, 4) = (4,2,2,1,1)$ so that

$\lambda =$

<table>
<tr><td> </td><td> </td><td> </td><td> </td></tr>
<tr><td> </td><td> </td></tr>
<tr><td> </td><td> </td></tr>
<tr><td> </td></tr>
<tr><td> </td></tr>
</table>

.

Then $\lambda \vdash 10$ is of $(2,s)$-hook shape since $\lambda_{2+1} = 2 \leq s$. The standard tableau $T \in \mathcal{HST}_{r,s}(\mathcal{L})$, where $\mathcal{L} = \{1,2,4,6,7,8,9,10,11,14\}$, depicted below

$T =$

<table>
<tr><td>1</td><td>6</td><td>11</td><td>14</td></tr>
<tr><td>2</td><td>8</td></tr>
<tr><td>4</td><td>9</td></tr>
<tr><td>7</td></tr>
<tr><td>10</td></tr>
</table>

determines two subtableaux

$T^{(1)} =$

<table>
<tr><td>1</td><td>6</td><td>11</td><td>14</td></tr>
<tr><td>2</td><td>8</td></tr>
</table>

and $T^{(2)} =$

<table>
<tr><td>4</td><td>7</td><td>10</td></tr>
<tr><td>9</td></tr>
</table>

of shapes $(4,2)$ and $(3,1) = (2,1,1)'$ respectively.

We have seen in Section 2 that the Brauer algebra $B_k(n-m)$ acts on $V^{\otimes k}$, and its action commutes with that of $\mathfrak{spo}(V,\beta)$ since it acts by $\mathfrak{spo}(V,\beta)$-module morphisms. In particular, the diagrams $d$ having all vertical edges form a subalgebra of $B_k(n-m)$ isomorphic to the group algebra of $S_k$, and the resulting action of $\mathbb{C}S_k$ on $V^{\otimes k}$ is by graded place permutations. We let $c_{p,q}$ denote the transformation on $V^{\otimes k}$ determined by the diagram in $B_k(n-m)$ with a horizontal edge connecting the $p$th and $q$th nodes on both the top and bottom, and with every other top node connected to the one directly below it. Thus, $c_{p,p+1}$ corresponds to the element $e_p$ in Section 2. We refer to the transformation $c_{p,q}$ as a *contraction mapping.* Clearly $c_{p,q} = c_{q,p}$, and we may suppose in working with contraction mappings that $p < q$. If $\underline{p} = \{p_1, \ldots, p_j\}$ and $\underline{q} = \{q_1, \ldots, q_j\}$ are disjoint ordered subsets of $\mathcal{K} = \{1, \ldots, k\}$, we set

$$c_{\underline{p},\underline{q}} = c_{p_1,q_1} \cdots c_{p_j,q_j} \quad \text{for } j = 1, \ldots, \lfloor k/2 \rfloor, \tag{3.5}$$

and we assume $c_{\emptyset,\emptyset}$ is the identity diagram. Since the image $\mathbb{C}\sum_{v,w \in B} F^{-1}_{v,w} v \otimes w$ of $e$ (see (2.17)) on $V^{\otimes 2}$ is a one-dimensional trivial $\mathfrak{spo}(V,\beta)$-module, we can draw the following conclusion (in Proposition 3.6 (a) below) that explains why these mappings are termed contractions:

**Proposition 3.6.** *(a)* $V^{\otimes k}c_{p,q} \cong V^{\otimes(k-2)}$ *for all* $p, q = 1, \ldots, k$ *with* $p \neq q$*.*

*(b) Let* $\overline{V^{\otimes k}}$ *be the subspace of all vectors in* $V^{\otimes k}$ *which are annihilated by all contractions* $c_{p,q}$*. Then* $\overline{V^{\otimes k}}$ *is an* $\mathfrak{spo}(V,\beta)$*-submodule and an* $S_k$*-submodule of* $V^{\otimes k}$*.*

*(c) Let* $\overline{V^{\otimes k}c_{\underline{p},\underline{q}}}$ *be the subspace of all vectors in* $V^{\otimes k}c_{\underline{p},\underline{q}}$ *which are annihilated by all contractions* $c_{s,t}$ *with* $s, t \in (\underline{p} \cup \underline{q})^c = \mathcal{K} - (\underline{p} \cup \underline{q})$*. Then* $\overline{V^{\otimes k}c_{\underline{p},\underline{q}}}$ *is a* $\mathfrak{spo}(V,\beta)$*-submodule (isomorphic to* $\overline{V^{\otimes(k-2j)}}$*,* $j = |\underline{p}| = |\underline{q}|$*) which is invariant under the action of the symmetric group* $S_{(\underline{p}\cup\underline{q})^c}$*.*

*Proof.* Since each $c_{p,q}$ commutes with the action of $\mathfrak{spo}(V,\beta)$, it determines an $\mathfrak{spo}(V,\beta)$-module endomorphism on $V^{\otimes k}$, and $\overline{V^{\otimes k}}$, which is the intersection of the kernels of the maps $c_{p,q}$, is thus an $\mathfrak{spo}(V,\beta)$-submodule. Its invariance under the action of $S_k$ comes from the identity

$$\sigma^{-1} c_{p,q} \sigma = c_{\sigma(p),\sigma(q)} \tag{3.7}$$

which holds for all $\sigma \in S_k$ and which can be verified in $B_k(n-m)$ by just multiplying the corresponding diagrams . Part (c) follows by combining (a) and (b). ∎

*Construction and linear independence of maximal vectors*

If $\underline{p} = \{p_1, \ldots, p_j\}$ and $\underline{q} = \{q_1, \ldots, q_j\}$ are disjoint ordered subsets of $\mathcal{K} = \{1, \ldots, k\}$, let $(\underline{p},\underline{q}) = \{(p_1,q_1), \ldots, (p_j,q_j)\}$, and denote by $\mathfrak{p}(j)$ the set of all such $(\underline{p},\underline{q})$. Assume

$$\mathfrak{p} = \bigcup_{j=0}^{\lfloor \frac{k}{2} \rfloor} \mathfrak{p}(j). \tag{3.8}$$

Using the contractions $c_{\underline{p},\underline{q}}$ and the symmetrizers $y_T$, we now construct maximal vectors of $V^{\otimes k}$.

**Theorem 3.9.** *Suppose $\lambda = (\lambda_1, \ldots, \lambda_\ell, 0, \ldots)$ is an $(r,s)$-hook shape partition of some integer $l$ such that $k - l = 2j$ and $0 \le j \le \lfloor k/2 \rfloor - 1$. Let $\lambda = \lambda_1\epsilon_1 + \cdots + \lambda_r\epsilon_r + \lambda'_{r+1}\delta_1 + \cdots + \lambda'_{r+s}\delta_s$ denote the weight determined by $\lambda$. Assume $(\underline{p}, \underline{q}) \in \mathfrak{p}(j)$ and fix $T \in \mathcal{HST}_{r,s}((\underline{p} \cup \underline{q})^c)$ of shape $\lambda$. Let $T^{(1)}$ and $T^{(2)}$ be the corresponding subtableaux of shapes $\lambda^{(1)} = (\lambda_1, \ldots, \lambda_r)$ and $\lambda^{(2)} = (\lambda_{r+1}, \ldots, \lambda_\ell)' = (\lambda'_{r+1}, \ldots, \lambda'_{r+s})$ respectively. Then $\theta = w_{T,\underline{p},\underline{q}} c_{\underline{p},\underline{q}} y_T$ is a maximal vector in $\overline{V^{\otimes k} c_{\underline{p},\underline{q}}} y_T$ of weight $\lambda$ where $w_{T,\underline{p},\underline{q}} = w_1 \cdots w_k$ is the simple tensor defined by*

$$w_i = \begin{cases} t_1 & \text{if } i \in \underline{p} \\ t_1^* & \text{if } i \in \underline{q} \\ t_j & \text{if } i \in (\underline{p} \cup \underline{q})^c \text{ and } i \text{ is in } j\text{th row of } T^{(1)} \\ u_j & \text{if } i \in (\underline{p} \cup \underline{q})^c \text{ and } i \text{ is in } j\text{th row of } T^{(2)}. \end{cases}$$

*If $n = 2s$ and $\ell(\lambda^{(2)}) = s$, then $\theta^\circ = w^\circ_{T,\underline{p},\underline{q}} c_{\underline{p},\underline{q}} y_T$ is a maximal vector of weight $\lambda^\circ = \lambda - 2\lambda'_{r+s}\delta_s$ where $w^\circ_{T,\underline{p},\underline{q}}$ is the simple tensor obtained from $w_{T,\underline{p},\underline{q}}$ by replacing the factors $u_s$ with $u_s^*$.*

*Proof.* Let $w = w_{T,\underline{p},\underline{q}}$. First we observe that $wc_{\underline{p},\underline{q}}$ is annihilated by all contractions having both subscripts in $(\underline{p} \cup \underline{q})^c$ because $\langle v, v' \rangle = 0$ whenever $v, v' \in \{t_1, \ldots, t_r, u_1, \ldots, u_s, (u_{s+1})\}$. Thus, $wc_{\underline{p},\underline{q}} \in \overline{V^{\otimes k} c_{\underline{p},\underline{q}}}$ and $\theta = wc_{\underline{p},\underline{q}} y_T \in \overline{V^{\otimes k} c_{\underline{p},\underline{q}}} y_T$. We also observe that the weight of $\theta$ is the same as the weight of $w$, which is $\lambda$, since the weights of $v$ and $v^*$ sum to zero. Now $w sgn(\psi)\psi\gamma = w$ for $\psi \in R_T$ and $\gamma \in C_T$ if and only if $\psi \in R_{T^{(1)}}$ and $\gamma \in R_{T^{(2)}}$. That assertion follows from the fact that $\beta(t_i, t_i) = 1$ for all $t_i$ and $\beta(u_i, u_i) = -1$ for all $u_i$, and the symmetric property of the bicharacter. Therefore, when $\theta$ is expressed as a linear combination of the basis tensors, the coefficient of $w$ in $\theta$ equals $h(\lambda)^{-1}|R_{T^{(1)}}| \cdot |R_{T^{(2)}}|$. Thus, $\theta \neq 0$. To prove that $\theta$ is a maximal vector, it suffices to show that it is annihilated by each of root vectors $x_j$ in (1.12). Since the image of a contraction is a trivial $\mathfrak{spo}(V,\beta)$-module, the sum of the terms with $x_j$ acting on a pair of contracted slots is zero. Thus, $x_j$ annihilates $wc_{\underline{p},\underline{q}}$ or produces a sum of tensors which are obtained by applying $x_j$ to a noncontracted factor. Each of those tensors has one factor in $(\underline{p} \cup \underline{q})^c$ whose subscript has been lowered by one or a factor $u_1$ which has been changed to $t_r$. For such tensors, each of them denoted by $v$, we argue that $vy_T = 0$.

First consider the case that $t_i$ has been changed to $t_{i-1}$ for $1 \le i \le r$ or that $u_1$ has been changed to $t_r$. Suppose $\psi \in R_T$. Then there exists $(a\ b) \in C_T$ which permutes two factors of $v\psi$ which are the same $t_j$ such that $v\psi(a\ b) = -v\psi$ (see (2.17)). Then $\sum_{\gamma \in C_T} sgn(\psi)(v\psi)\gamma = \sum_{\gamma \in C_T} sgn(\psi)(v\psi)(a\,b)\gamma = -\sum_{\gamma \in C_T} sgn(\psi)(v\psi)\gamma$. Thus, we have $\sum_{\gamma \in C_T} sgn(\psi)(v\psi)\gamma = 0$ for each $\psi \in R_T$ and $vy_T = 0$.

Now consider the case that $u_i$ has been changed to $u_{i-1}$. Then there exists $(a\ b) \in R_T$ permuting two factors of $v$ which are the same $u_i$ such that $v(a\ b) = v$ (see (2.17)). As a result, $\sum_{\psi \in R_T} sgn(\psi)v\psi = \sum_{\psi \in R_T} sgn((a\ b))sgn(\psi)(v(a\ b))\psi = -\sum_{\psi \in R_T} sgn(\psi)v\psi$. Therefore, $\sum_{\psi \in R_T} sgn(\psi)v\psi = 0$, and $vy_T = 0$. Thus, for all such $v$ we have $vy_T = 0$, and consequently $x_j(wc_{\underline{p},\underline{q}} y_T) = (x_j wc_{\underline{p},\underline{q}})y_T = 0$.

If $n = 2s$ and $\ell(\lambda^{(2)}) = s$, then the same argument as above proves that $\theta^\circ$ is a maximal vector of weight $\lambda - 2\lambda'_{r+s}\delta_s$ since in this case $x_j$ for all $j \neq r+s-1, r+s$ acts on $w^\circ$ exactly as it acted on $w$. The action of $x_{r+s-1}$ is zero since all the factors $u_s$ have been transformed to $u_s^*$, while the action of $x_{r+s}$ on $w^\circ$ produces a sum of tensors each of which has one $u_s^*$ factor changed to $u_{s-1}$. Consequently, $\theta^\circ$ is a maximal vector as claimed. ∎

Now consider $V^{\otimes k} c_{\underline{p},\underline{q}}$ where $(\underline{p},\underline{q}) \in \mathfrak{p}(\lfloor k/2 \rfloor)$. When $k$ is even, $V^{\otimes k} c_{\underline{p},\underline{q}}$ is the trivial $\mathfrak{spo}(V,\beta)$-module $\kappa = \mathbf{1}$, and when $k$ is odd, it is isomorphic to $V$. If $k$ is odd, then $T \in \mathcal{HST}_{r,s}((\underline{p} \cup \underline{q})^c)$ consists of a single box with a number in it so that $y_T = id$. When $k$ is even, then $(\underline{p} \cup \underline{q})^c = \emptyset$, $T = \emptyset$, and we set $y_\emptyset = id$. Let $w_{T,\underline{p},\underline{q}} = w_1 \cdots w_k$ be the simple tensor prescribed by Theorem 3.9 in these cases. Then $w_{T,\underline{p},\underline{q}} c_{\underline{p},\underline{q}} y_T$ is a maximal vector of weight 0 when $k$ is even, and $w_{T,\underline{p},\underline{q}} c_{\underline{p},\underline{q}} y_T$ is a maximal vector of weight $\epsilon_1$ when $k$ is odd. Moreover, $V^{\otimes k} c_{\underline{p},\underline{q}}$ is generated by $w_{T,\underline{p},\underline{q}} c_{\underline{p},\underline{q}} y_T$. To summarize we have:

**Corollary 3.10.** *For $(\underline{p},\underline{q}) \in \mathfrak{p}(\lfloor k/2 \rfloor)$*

$$V^{\otimes k} c_{\underline{p},\underline{q}} \cong \begin{cases} \kappa & \text{if } k \text{ is even} \\ V & \text{if } k \text{ is odd,} \end{cases}$$

*and $w_{T,\underline{p},\underline{q}} c_{\underline{p},\underline{q}} y_T$ is a maximal vector of $V^{\otimes k} c_{\underline{p},\underline{q}}$ of weight 0 when $k$ is even and of weight $\epsilon_1$ when $k$ is odd, where $w_{T,\underline{p},\underline{q}}$ is as in Theorem 3.9.*

Now we show that the maximal vectors constructed above are linearly independent if the rank of the algebra $\mathfrak{spo}(V,\beta)$ is sufficiently large. Suppose $r+s \geq k$ and $s \geq 2$. Assume $\lambda = (\lambda_1, \ldots, \lambda_\ell, 0, \ldots)$ is an $(r,s)$-hook shape partition of $l$ for some $0 \leq l \leq k$ and let $\lambda^{(1)}$ and $\lambda^{(2)}$ be the associated subpartitions. We argue that in this case $\ell(\lambda^{(2)}) < s$. Indeed, if $\ell(\lambda^{(2)}) = s$, then since $\lambda_r \geq \ell(\lambda^{(2)}) = s$, $|\lambda| = |\lambda^{(1)}| + |\lambda^{(2)}| \geq rs + s > r + s \geq k$. This contradicts the fact that $\lambda$ partitions $l \leq k$. Thus, when $r+s \geq k$, we need only work with the maximal vectors of the form $w_{T,\underline{p},\underline{q}} c_{\underline{p},\underline{q}} y_T$.

**Theorem 3.11.** *Suppose that $r+s \geq k$. Then*

$$\{w_{T,\underline{p},\underline{q}} c_{\underline{p},\underline{q}} y_T \mid (\underline{p},\underline{q}) \in \mathfrak{p},\ T \in \mathcal{HST}_{r,s}((\underline{p} \cup \underline{q})^c)\}$$

*is a linearly independent set of maximal vectors.*

*Proof.* Vectors corresponding to different numbers of contractions are linearly independent because those vectors have different weights. As the weight is determined by the shape of the tableau $T$, vectors corresponding to tableaux with different shapes are linearly independent. Thus, it suffices to consider the $(\underline{p},\underline{q})$'s in $\mathfrak{p}(j)$ for some fixed value of $j$, and the tableau $T$ of a fixed shape $\lambda$. Suppose

$$\sum_{(\underline{p},\underline{q}) \in \mathfrak{p}(j)} \sum_{T \in \mathcal{HST}_{r,s}((\underline{p} \cup \underline{q})^c)} a_{T,\underline{p},\underline{q}} w_{T,\underline{p},\underline{q}} c_{\underline{p},\underline{q}} y_T = 0,$$

where the sum is taken over all the tableaux $T \in \mathcal{HST}_{r,s}((\underline{p} \cup \underline{q})^c)$ of some given shape $\lambda = (\lambda_1, \ldots, \lambda_\ell, 0, \ldots)$, which is an $(r,s)$-hook shape partition. We suppose that $\lambda^{(1)} = (\lambda_1, \ldots, \lambda_r)$ and $\lambda^{(2)} = (\lambda_{r+1}, \ldots, \lambda_\ell)' = (\lambda'_{r+1}, \ldots, \lambda'_{r+s})$ are the partitions determined by $\lambda$, and let $\ell_1 = \ell(\lambda^{(1)}$ and $\ell_2 = \ell(\lambda^{(2)})$. Among the simple tensors involved in the expression for a fixed $w_{T,\underline{p},\underline{q}} c_{\underline{p},\underline{q}}$ there is a unique tensor $\zeta_{T,\underline{p},\underline{q}}$ matching the following description:

First when $\ell_2 = 0$, the sequence of positions $p_1, q_1, p_2, q_2, \ldots, p_j, q_j$ in $\zeta_{T,\underline{p},\underline{q}}$ contains the first $2j$ terms of the sequence $t_{\ell_1+1}, t^*_{\ell_1+1}, t_{\ell_1+2}, t^*_{\ell_1+2}, \ldots, t_r, t^*_r, u_1, u^*_1, \ldots, u_s, u^*_s$. Note there are sufficiently many vectors to make this work since $r+s-j \geq k-j \geq k-2j = |\lambda| \geq \ell_1$ implies $r+s-\ell_1 \geq j$. Similarly when $\ell_2 \neq 0$, then the sequence of positions $p_1, q_1, p_2, q_2, \ldots, p_j, q_j$ in $\zeta_{T,\underline{p},\underline{q}}$ contains the first $2j$ terms of the sequence $u_{\ell_2+1}, u^*_{\ell_2+1}, u_{\ell_2+2}, u^*_{\ell_2+2}, \ldots, u_s, u^*_s$, and there are enough vectors in this case since $r+s-j \geq k-j \geq k-2j = |\lambda| \geq \ell_1 + \ell_2 = r + \ell_2$ so that $s - \ell_2 \geq j$.

Then it follows that

$$\sum_{(\underline{p},\underline{q})\in\mathfrak{p}(j)} \sum_{T\in\mathcal{HST}_{r,s}((\underline{p}\cup\underline{q})^c)} a_{T,\underline{p},\underline{q}}\zeta_{T,\underline{p},\underline{q}}y_T = 0,$$

where the sum is over $T$ of shape $\lambda$. The simple tensors in the above corresponding to different choices of $(\underline{p},\underline{q})$ are linearly independent, for the positions of the factors $t_{\ell(\lambda^{(1)})+i}$ and $t^*_{\ell(\lambda^{(1)})+i}$, or $u_{\ell(\lambda^{(2)})+i}$ and $u^*_{\ell(\lambda^{(2)})+i}$ are sufficient to distinguish them. Thus for any $(\underline{p},\underline{q}) \in \mathfrak{p}(j)$,

$$\sum_{T\in\mathcal{HST}_{r,s}((\underline{p}\cup\underline{q})^c)} a_{T,\underline{p},\underline{q}}\zeta_{T,\underline{p},\underline{q}}y_T = 0. \tag{3.12}$$

If some coefficient in (3.12) is nonzero, then since the sum is over tableaux of the same shape, we may choose $T'$ minimal in the order with the property that $a_{T',\underline{p},\underline{q}} \neq 0$. By Lemma 3.3, $T > T'$ implies that $y_T y_{T'} = 0$. Thus applying $y_{T'}$ to (3.12), we obtain

$$0 = a_{T',\underline{p},\underline{q}}\zeta_{T',\underline{p},\underline{q}}y_{T'}y_{T'} = a_{T',\underline{p},\underline{q}}\zeta_{T',\underline{p},\underline{q}}y_{T'}.$$

Since the vector $\zeta_{T',\underline{p},\underline{q}}$ appears in $\zeta_{T',\underline{p},\underline{q}}y_{T'}$ with coefficient equal to $h(\lambda)^{-1}|R_{T'^{(1)}}|\cdot|R_{T'^{(2)}}|$ as we discussed in the proof of Theorem 3.9, $\zeta_{T',\underline{p},\underline{q}}y_{T'} \neq 0$. Hence, no such $T'$ exists, and $a_{T,\underline{p},\underline{q}} = 0$ for all $T$ and all $(\underline{p},\underline{q})$. ∎

## 4. The $\mathfrak{spo}(V,\beta)$-modules $T^\lambda$ and their characters

*Definition of the modules $T^\lambda$*

Let $\mathfrak{spo}(V,\beta)$ be the orthosymplectic Lie color algebra and assume that $m = \dim V_{(0)}$ and $n = \dim V_{(1)}$. Recall that by Theorem 2.16 there is an action (which we write on the right) of the Brauer algebra $B_k(n-m)$ on $V^{\otimes k}$ which commutes with the action of $\mathfrak{spo}(V,\beta)$ on $V^{\otimes k}$. For convenience we write $w\,d$ rather than $w\,\Psi_d$ for the action of the element $d \in B_k(n-m)$ on $w \in V^{\otimes k}$.

Suppose that $|n-m| > k$. Then by the work of Wenzl [W] we know that the Brauer algebra $B_k(n-m)$ is semisimple and has simple summands indexed by the partitions in the set

$$\widehat{B}_k = \{\lambda \vdash k-2h \mid h = 0, 1, \ldots, \lfloor k/2 \rfloor\}.$$

More specifically, there are positive integers $d_\lambda$, $\lambda \in \widehat{B}_k$, and an isomorphism

$$\Gamma\colon B_k(n-m) \longrightarrow \bigoplus_{\lambda\in\widehat{B}_k} M_{d_\lambda}(\kappa),$$

where $M_{d_\lambda}(\kappa)$ denotes the full matrix algebra of $d_\lambda \times d_\lambda$ matrices with entries in $\kappa$.

For each $\lambda \in \widehat{B}_k$ and $1 \le P, Q \le d_\lambda$, let $E^\lambda_{P,Q}$ denote the matrix unit in the $\lambda$th block of $\bigoplus_{\lambda \in \widehat{B}_k} M_{d_\lambda}(\kappa)$ which has a 1 as the $(P,Q)$ entry and zeros everywhere else. Let $e^\lambda_{P,Q} = \Gamma^{-1}(E^\lambda_{P,Q})$ and define

$$T^{\lambda'} = V^{\otimes k} e^\lambda_{P,Q}, \tag{4.1}$$

where $\lambda'$ denotes the partition conjugate to $\lambda$, and the action of $B_k(n-m)$ on $V^{\otimes k}$ is as given in Theorem 2.16. Then $T^{\lambda'}$ is a $\mathfrak{g}$-module since the action of $\mathfrak{g}$ on $V^{\otimes k}$ commutes with the action of $e^\lambda_{P,Q} \in B_k(n-m)$; however, the module $T^{\lambda'}$ may be $(0)$.

**Proposition 4.2.** *If $|n-m| > k$, then the following hold:*

*(a) The $\mathfrak{g}$-module $T^{\lambda'}$ is well-defined (up to $\mathfrak{g}$-module isomorphism it does not depend on the choice of $P$ and $Q$).*

*(b) As modules for $\mathfrak{g} \times B_k(n-m)$,*

$$V^{\otimes k} \cong \bigoplus_{\lambda \in \widehat{B}_k} T^{\lambda'} \otimes B_\lambda \cong \bigoplus_{\lambda \in \widehat{B}_k} T^{\lambda} \otimes B_{\lambda'}$$

*where $B_\lambda$ is the irreducible $B_k(n-m)$-module labeled by $\lambda$ for each $\lambda \in \widehat{B}_k$.*

*Proof.* We shall only sketch the proof of this standard result.

(a) The map

$$\begin{array}{llll} e^\lambda_{Q,S}\colon & V^{\otimes k} e^\lambda_{P,Q} & \to & V^{\otimes k} e^\lambda_{R,S} \\ & w\, e^\lambda_{P,Q} & \mapsto & w\, e^\lambda_{P,Q} e^\lambda_{Q,S} \quad = (w\, e^\lambda_{P,R}) e^\lambda_{R,S} \end{array}$$

is a $\mathfrak{g}$-module isomorphism with inverse given by

$$\begin{array}{llll} e^\lambda_{S,Q}\colon & V^{\otimes k} e^\lambda_{R,S} & \to & V^{\otimes k} e^\lambda_{P,Q} \\ & w\, e^\lambda_{R,S} & \mapsto & w\, e^\lambda_{R,S} e^\lambda_{S,Q} \quad = (w\, e^\lambda_{R,P}) e^\lambda_{P,Q}. \end{array}$$

(b) For each $\lambda \in \widehat{B}_k$, the irreducible $B_k(n-m)$-module $B_\lambda$ has a basis $\{e^\lambda_Q \mid 1 \le Q \le d_\lambda\}$, and the action of $B_k(n-m)$ on this space is given by

$$e^\lambda_Q e^\mu_{R,S} = \delta_{\lambda,\mu} \delta_{Q,R} e^\lambda_S.$$

Note that since $1 = \sum_{\lambda \in \widehat{B}_k} \sum_{1 \le Q \le d_\lambda} e^\lambda_{Q,Q}$ we may write

$$V^{\otimes k} = \sum_{\lambda \in \widehat{B}_k} \sum_{1 \le Q \le d_\lambda} V^{\otimes k} e^\lambda_{Q,Q}.$$

For convenience let us fix a particular index, say $P = 1$, for each $\lambda \in \widehat{B}_k$ and set $T^{\lambda'} = V^{\otimes k} e^\lambda_{1,1}$. Now one can check that the maps

$$\begin{array}{llll} \Phi\colon & \bigoplus_\lambda T^{\lambda'} \otimes B_\lambda & \to & V^{\otimes k} = \sum_\lambda \sum_Q V^{\otimes k} e^\lambda_{Q,Q} \\ & w\, e^\lambda_{1,1} \otimes e^\lambda_Q & \mapsto & w\, e^\lambda_{1,Q} = (w\, e^\lambda_{1,Q}) e^\lambda_{Q,Q} \quad \in V^{\otimes k} e^\lambda_{Q,Q}, \end{array}$$

and

$$\begin{array}{rlcl}\Phi'\colon & V^{\otimes k}=\sum_\lambda\sum_Q V^{\otimes k}e^\lambda_{Q,Q} & \to & \bigoplus_\lambda T^{\lambda'}\otimes B_\lambda \\ & w=\sum_\lambda\sum_Q w\,e^\lambda_{Q,Q} & \mapsto & \sum_\lambda\sum_Q (w\,e^\lambda_{Q,1})e^\lambda_{1,1}\otimes e^\lambda_Q\end{array}$$

are $\mathfrak{g}\times B_k(n-m)$-module homomorphisms which are inverses of each other. ∎

*Characters of the Brauer algebras*

If $d_1$ is an $m_1$-diagram and $d_2$ is an $m_2$-diagram, then $d_1\otimes d_2$ is the $(m_1+m_2)$-diagram obtained by placing $d_1$ to the right of $d_2$. Let $e$ denote the 2-diagram

$e\ =$ 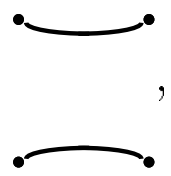 , (4.3)

and let $\gamma_m$ denote the $m$-diagram

$\gamma_m\ =$ 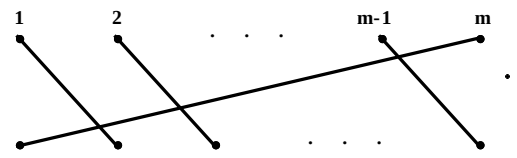

.

For a partition $\mu=(\mu_1,\mu_2,\ldots,\mu_t)$, let $\gamma_\mu=\gamma_{\mu_1}\otimes\gamma_{\mu_2}\otimes\cdots\otimes\gamma_{\mu_t}$. We have the following results from [R].

**Theorem 4.4.**
*(a) If $\eta\neq 0$, then any character of the Brauer algebra $B_k(\eta)$ is completely determined by its values on the elements $e^{\otimes j}\otimes\gamma_\mu$, for $j=0,1,\ldots,\lfloor k/2\rfloor$ and $\mu\vdash k-2j$.*
*(b) Assume that $\kappa\subseteq\mathbb{C}$ so that $\eta\in\mathbb{C}$. If $\eta\notin\mathbb{Z}$ or if $\eta\in\mathbb{Z}$ and $|\eta|>k$, then the irreducible characters of $B_k(\eta)$,*

$$\chi^\lambda_{B_k(\eta)}\colon B_k(\eta)\to\kappa,$$

*which are indexed by the partitions $\lambda$ in the set $\widehat{B}_k=\{\lambda\vdash k-2h\mid h=0,1,\ldots,\lfloor k/2\rfloor\}$, are given by*

$$\chi^\lambda_{B_k(\eta)}(e^{\otimes j}\otimes\gamma_\mu)=\eta^j\sum_{\nu\supseteq\lambda}\left(\sum_{\theta\text{ even}}c^\nu_{\lambda,\theta}\right)\chi^\nu_{S_{k-2j}}(\mu),$$

*where $\chi^\nu_{S_{k-2j}}(\mu)$ denotes the irreducible character of the symmetric group $S_{k-2j}$ labeled by $\nu\vdash k-2j$ evaluated at the conjugacy class indexed by $\mu\vdash k-2j$, and $c^\nu_{\lambda,\theta}$ is the Littlewood-Richardson coefficient (see (4.9) below), and the inner sum is over all partitions $\theta$ with even parts.*

*Weighted traces*

Assume that the basis $B=B_0\cup B_1$ of $V$ is as in (1.8). Let $\{z_b\}_{b\in B}$ be a set of commuting variables indexed by the elements of $B$ such that

$$z_{b^*}=z_b^{-1},\qquad\text{for all } b\in B,\tag{4.5}$$

(in particular, $z_{u_{s+1}} = 1$). Define an operator $D$ on $V^{\otimes k}$ by

$$D(b_1 \otimes b_2 \otimes \cdots \otimes b_k) = z_{b_1} z_{b_2} \cdots z_{b_k} (b_1 \otimes b_2 \otimes \cdots \otimes b_k),$$

for $b_1, b_2, \dots b_k \in B$, or simply $D(b_1 \cdots b_k) = z_{b_1} z_{b_2} \cdots z_{b_k}(b_1 \cdots b_k)$ deleting the tensor symbols. Since (by a brute force check) $D$ commutes with the action of the generators $e_i$ and $s_i$, $1 \leq i \leq k-1$ of the Brauer algebra, it follows that $D$ commutes with the action of the Brauer algebra $B_k(n-m)$ on $V^{\otimes k}$. Using $D$ we define the *weighted trace* of $d \in B_k(n-m)$ by

$$\mathrm{wtr}(d) = \mathrm{Tr}_{V^{\otimes k}}(dD). \tag{4.6}$$

**Lemma 4.7.**

(a) *If $d$ is a $k$-diagram, then*

$$\mathrm{wtr}(d) = \sum_{b_1, b_2, \dots, b_k \in B} z_{b_1} z_{b_2} \cdots z_{b_k} d_{\underline{b}, \underline{b}},$$

*where $\underline{b} = (b_1, b_2, \dots, b_k)$ and $d_{\underline{b},\underline{b}}$ denotes the weight of the diagram $d$ labeled on the top by $b_1, \dots, b_k \in B$ and on the bottom also by $b_1, \dots, b_k$.*

(b) *The map* wtr *is a trace on $B_k(n-m)$, i.e.,* $\mathrm{wtr}(d_1 d_2) = \mathrm{wtr}(d_2 d_1)$, *for all $d_1, d_2 \in B_k(n-m)$.*

(c) *If $d_1$ is a $k$-diagram and $d_2$ is an $\ell$-diagram then* $\mathrm{wtr}(d_1 \otimes d_2) = \mathrm{wtr}(d_1)\mathrm{wtr}(d_2)$.

(d) *If $d \in B_k(n-m)$ then*

$$\mathrm{wtr}(d) = \sum_{\lambda \in \widehat{B}_k} \chi^{\lambda'}_{B_k(n-m)}(d) \mathrm{char}(T^\lambda),$$

*where $\chi^{\lambda'}_{B_k(n-m)}$ denotes the irreducible character of $B_k(n-m)$ indexed by $\lambda' \in \widehat{B}_k$ and* $\mathrm{char}(T^\lambda) = \mathrm{Tr}_{T^\lambda}(D)$.

*Proof.* (a) This is a direct consequence of the definition of the weighted trace.

(b) Since the actions of $D$ and $B_k(n-m)$ commute on $V^{\otimes k}$, it follows that

$$\mathrm{wtr}(d_1 d_2) = \mathrm{Tr}(d_1 d_2 D) = \mathrm{Tr}(d_2 D d_1) = \mathrm{Tr}(d_2 d_1 D) = \mathrm{wtr}(d_2 d_1).$$

(c) This comes from (a) and the fact that if $b_1, \dots, b_k, b_{k+1}, \dots, b_{k+\ell} \in B$, then

$$D(b_1 \cdots b_k b_{k+1} \cdots b_\ell) = \Big(D(b_1 \cdots b_k)\Big)\Big(D(b_{k+1} \cdots b_{k+\ell})\Big)$$

and

$$\Big(b_1 \cdots b_k b_{k+1} \cdots b_{k+\ell}\Big)(d_1 \otimes d_2) = \Big((b_1 \cdots b_k) d_1\Big)\Big((b_{k+1} \cdots b_{k+\ell}) d_2\Big).$$

(d) It is clear that $T^\lambda$ is also invariant under the action of $D$ as the actions of $D$ and the Brauer algebra $B_k(n-m)$ commute. Thus, all the statements of Proposition 4.2(b) hold with $D$ in place of $\mathfrak{g}$. By taking the trace of the operator $dD$ on each side of the isomorphism in Proposition 4.2(b) (with $\mathfrak{g}$ replaced by $D$) we have that

$$\mathrm{wtr}(d) = \mathrm{Tr}_{V^{\otimes k}}(dD) = \sum_{\lambda\in\widehat{B}_k} \mathrm{Tr}_{B_{\lambda'}}(d)\mathrm{Tr}_{T^\lambda}(D) = \sum_{\lambda\in\widehat{B}_k} \chi^{\lambda'}_{B_k(n-m)}(d)\mathrm{char}(T^\lambda). \qquad \blacksquare$$

In view of Lemma 4.7(c) and Theorem 4.4(b), it is sufficient to compute the weighted traces $\mathrm{wtr}(d)$ when $d$ is a $k$-diagram of the form $e^{\otimes j}\otimes\gamma_\mu$. This is done in the following lemma.

**Lemma 4.8.** *Let $\mu=(\mu_1,\mu_2,\ldots,\mu_t)$ be a partition of $k-2j$ for some $j\in\{0,1,\ldots,\lfloor k/2\rfloor\}$. Let $e^{\otimes j}\otimes\gamma_\mu$ be the $k$-diagram coming from (4.3). For each positive integer $\ell$ suppose*

$$\overline{p}_\ell(\tilde{Z}) = \sum_{b\in B_1} z_b^\ell - \sum_{b\in B_0}(-1)^\ell z_b^\ell, \quad \text{and let} \quad \overline{p}_\mu(\tilde{Z}) = \overline{p}_{\mu_1}(\tilde{Z})\cdots\overline{p}_{\mu_t}(\tilde{Z}).$$

*Then*

$$\mathrm{wtr}(e^{\otimes j}\otimes\gamma_\mu) = (n-m)^j\overline{p}_\mu(\tilde{Z}).$$

*Proof.* From the definition of the labeled diagrams and the properties of the matrices $F_B$ and $F_B^{-1}$ in (1.9a,b) we see that

$$e_{\underline{b},\underline{b}} = F^{-1}_{b_1,b_2}F_{b_1,b_2} = \begin{cases} -1 & \text{if } b_1^*=b_2\in B_0,\\ 1 & \text{if } b_1^*=b_2\in B_1,\\ 0 & \text{otherwise,}\end{cases}$$

where $\underline{b}=(b_1,b_2)$. If $\underline{b}=(b_1,b_2,\ldots,b_k)$ then

$$(\gamma_\ell)_{\underline{b},\underline{b}} = \prod_{i=1}^{\ell-1}\delta_{b_i,b_{i+1}}\big(-\beta(b_\ell,b_i)\big) = \begin{cases}(-1)^{\ell-1}\beta(b,b)^{\ell-1} & \text{if } \underline{b}=(b,b,\ldots,b),\\ 0 & \text{otherwise.}\end{cases}$$

Consequently, we have

$$\begin{aligned}\mathrm{wtr}(e) &= \sum_{b_1,b_2\in B} z_{b_1}z_{b_2}e_{\underline{b},\underline{b}} = \sum_{b\in B_0}(-1)z_bz_{b^*} + \sum_{b\in B_1} z_bz_{b^*} = \sum_{b\in B_0}(-1) + \sum_{b\in B_1} 1 = n-m,\\ \mathrm{wtr}(\gamma_\ell) &= \sum_{b\in B_0} z_b^\ell(-1)^{\ell-1} + \sum_{b\in B_1} z_b^\ell.\end{aligned}$$

The result now follows since by Lemma 4.7(c) we have that

$$\mathrm{wtr}(e^{\otimes j}\otimes\gamma_\mu) = \mathrm{wtr}(e)^j\mathrm{wtr}(\gamma_{\mu_1})\mathrm{wtr}(\gamma_{\mu_2})\cdots\mathrm{wtr}(\gamma_{\mu_t}). \qquad \blacksquare$$

*Symmetric functions and a combinatorial description of* $\mathrm{char}(T^\lambda)$

We adopt the notation for partitions and symmetric functions found in [M]. Let us begin by recalling a few definitions.

Suppose $Y = \{y_1, \ldots, y_q\}$ is a set of commuting variables ordered by $y_1 < y_2 < \ldots < y_q$, and assume $\lambda$ and $\mu$ are partitions such that $\mu \subseteq \lambda$. A *column-strict tableau* of shape $\lambda/\mu$ is a filling of the boxes in the Ferrers diagram of $\lambda/\mu$ with $y_j$'s such that the $y_j$'s are weakly increasing (left to right) across the rows and are strictly increasing down the columns. Associated to a column-strict tableau $T$ of shape $\lambda/\mu$ is the monomial $y^T$, which is the product over all boxes of $\lambda/\mu$ of the elements $y_j$ in the boxes. The *skew Schur function* is the sum

$$s_{\lambda/\mu}(Y) = \sum_T y^T,$$

over all column-strict tableaux $T$ of shape $\lambda/\mu$. The (ordinary) Schur function $s_\lambda(Y)$ is just the skew Schur function $s_{\lambda/\mu}(Y)$ in the special case that $\mu = \emptyset$. The *Littlewood-Richardson coefficient* is the nonnegative integer $c^\lambda_{\mu,\nu}$ defined by

$$s_\mu(Y)s_\nu(Y) = \sum_\lambda c^\lambda_{\mu,\nu} s_\lambda(Y). \tag{4.9}$$

The Littlewood-Richardson coefficients exhibit certain symmetries such as

$$c^\lambda_{\mu,\nu} = c^{\lambda'}_{\mu',\nu'}. \tag{4.10}$$

It is a standard fact that the skew Schur functions can be written in terms of the ordinary Schur functions via the relation

$$s_{\lambda/\mu}(Y) = \sum_\nu c^\lambda_{\mu,\nu} s_\nu(Y). \tag{4.11}$$

A partition $\rho$ can be specified by its Frobenius notation $\rho = (q_1, \ldots, q_p | r_1, \ldots, r_p)$ where the main diagonal of $\rho$ consists of $p$ boxes $(i,i)$, $1 \le i \le p$, and there are $q_i$ boxes to the right of $(i,i)$ in the $i$th row and $r_i$ boxes below $(i,i)$ in the $i$th column. When $\rho$ partitions $k$, it is customary to write $|\rho| = k$. We have the following identities involving Schur functions (see [M, pp. 77-79]):

$$\prod_{1\le i\le j\le q} (1 - y_i y_j) = \sum_\rho (-1)^{|\rho|/2} s_\rho(Y), \tag{4.12}$$

where the sum is over all partitions of the form $\rho = (r_1+1, \ldots, r_p+1 \mid r_1, \ldots, r_p)$, $r_1 \le q-1$, in Frobenius notation;

$$\prod_{1\le i<j\le q} (1 - y_i y_j) = \sum_\pi (-1)^{|\pi|/2} s_\pi(Y), \tag{4.13}$$

summed over the partitions $\pi = (r_1 - 1, \ldots, r_p - 1 \mid r_1, \ldots, r_p)$, $r_1 \le q-1$, in Frobenius notation; and

$$\prod_{1\le i\le j\le q} (1 - y_i y_j)^{-1} = \sum_{\theta \text{ even}} s_\theta(Y), \tag{4.14}$$

where the sum is over all partitions $\theta$ with even length rows (all parts even). By using the above identities one can easily prove the following

**Proposition 4.15.**

$$\sum_{\nu}\left(\sum_{\rho}(-1)^{|\rho|/2}c^{\tau}_{\nu,\rho}\right)\left(\sum_{\theta\text{ even}}c^{\nu}_{\lambda,\theta}\right)=\delta_{\tau,\lambda},\tag{4.16}$$

*where the sums over* $\rho$ *and* $\theta$ *are as in (4.12) and (4.14), respectively.*

*Proof.* By using (4.12), (4.14), and (4.9) to expand, we get

$$\begin{aligned}s_\lambda(Y)&=s_\lambda(Y)\prod_{1\le i\le j\le q}(1-y_iy_j)^{-1}\prod_{1\le i\le j\le q}(1-y_iy_j)\\&=\sum_{\tau}s_\tau(Y)\sum_{\nu}\left(\sum_{\theta\text{ even}}c^{\nu}_{\lambda,\theta}\right)\left(\sum_{\rho}(-1)^{|\rho|/2}c^{\tau}_{\nu,\rho}\right)\end{aligned}$$

and by comparing coefficients of the Schur functions on each sides, we get identity (4.16). ∎

Let us recall the definition of the *hook Schur functions* given in [BR] and [Se], which describe the characters of the irreducible $\mathrm{gl}(V,\beta)$-modules appearing in $V^{\otimes k}$ when $G=\mathbb{Z}_2$, and $\beta(a,b)=(-1)^{ab}$; that is, for the general linear superalgebra. Order the variables $z_b$, $b\in B=B_0\cup B_1$, by

$$z_{t_1}<z_{t_1^*}<\cdots<z_{t_r}<z_{t_r^*}<z_{u_1}<z_{u_1^*}<\cdots<z_{u_s}<z_{u_s^*}<(z_{u_{s+1}}).$$

A *bitableau* of shape $\lambda$ is a filling of the Ferrers diagram of $\lambda$ with elements of $B$ such that

(1) the portion of the diagram filled with $z_t$'s is the diagram of a partition $\mu\subseteq\lambda$,
(2) the $z_t$'s are weakly increasing (left to right) along rows and strictly increasing down columns.
(3) the $z_u$'s are weakly increasing down the columns and strictly increasing left to right across the rows.

If $T$ is a bitableau of shape $\lambda$, then $z^T$ is the product over all boxes of $\lambda$ of the elements $z_b$ in the boxes, and the *hook Schur function* (see [M, Chap. I §3 Ex. 23-24 and Chap. I §5 Ex. 23]), is given by

$$s_\lambda(\tilde Z)=\sum_T z^T,$$

where the sum is over all bitableaux of shape $\lambda$.

The hook Schur functions satisfy the following identities (see [M, Chap. I (4.3), (4.3'), (7.7) and §3 Ex. 23]):

$$\sum_\lambda s_\lambda(\tilde Z)s_\lambda(Y)=\frac{\prod_{b\in B_1}\prod_{j=1}^q(1+z_by_j)}{\prod_{b\in B_0}\prod_{j=1}^q(1-z_by_j)};\tag{4.17}$$

$$\sum_\lambda s_{\lambda'}(\tilde Z)s_\lambda(Y)=\frac{\prod_{b\in B_0}\prod_{j=1}^q(1+z_by_j)}{\prod_{b\in B_1}\prod_{j=1}^q(1-z_by_j)},\tag{4.18}$$

where $s_\lambda(Y)$ denotes the ordinary Schur function in the variables $y_1, \ldots, y_q$. If $\lambda \vdash k$, then

$$s_\lambda(\tilde{Z}) = \sum_{\mu \vdash k} \frac{\chi^\lambda_{S_k}(\mu)}{\mathfrak{z}_\mu} p_\mu(\tilde{Z}), \tag{4.19}$$

where (i) $\chi^\lambda_{S_k}(\mu)$ denotes the irreducible character of the symmetric group indexed by the partition $\lambda \vdash k$ evaluated at the conjugacy class labeled by the partition $\mu \vdash k$, (ii) $\mathfrak{z}_\mu$ is the cardinality of the centralizer of an element of cycle type $\mu$ in $S_k$, and (iii) $p_\mu(\tilde{Z})$ is the function which is defined by

$$p_\ell(\tilde{Z}) = \sum_{b \in B_0} z_b^\ell - \sum_{b \in B_1} (-1)^\ell z_b^\ell, \quad \text{and} \quad p_\mu(\tilde{Z}) = p_{\mu_1}(\tilde{Z}) \cdots p_{\mu_t}(\tilde{Z}).$$

The functions $\overline{p}_\mu(\tilde{Z})$ defined in Lemma 4.8 satisfy $\overline{p}_\ell(\tilde{Z}) = (-1)^\ell(-p_\ell(\tilde{Z}))$ for all $\ell > 0$, and $\overline{p}_\mu(\tilde{Z}) = \prod_i (-1)^{\mu_i}(-p_{\mu_i}(\tilde{Z})) = (-1)^{k-\ell(\mu)} p_\mu(\tilde{Z})$, for a partition $\mu \vdash k$. Then (4.19) and the fact that $\chi^\lambda_{S_k}(\mu) = (-1)^{k-\ell(\mu)} \chi^{\lambda'}_{S_k}(\mu)$, (see for example, [M, Chap. 1, §7, Ex. 2]) gives

$$s_\lambda(\tilde{Z}) = \sum_{\mu \vdash k} \frac{\chi^{\lambda'}_{S_k}(\mu)}{\mathfrak{z}_\mu} \overline{p}_\mu(\tilde{Z}). \tag{4.20}$$

Let $Z_0 = \{z_b \mid b \in B_0\}$ be ordered as in the definition of the hook Schur function. A *symplectic tableau of shape* $\lambda$ is a filling of the boxes in the Ferrers diagram of $\lambda$ with $z_b$'s, $b \in B_0$, such that

(1) the $z_b$'s are weakly increasing (left to right) along rows and strictly increasing down columns,

(2) the elements $z_{t_i}$ and $z_{t_i^*}$ never appear in a row with number greater than $i$.

Associated to a symplectic tableau $T$ of shape $\lambda$ is the product $z^T$ of all the elements $z_b$ in the boxes of $\lambda$. The *symplectic Schur function* is defined by

$$sp_\lambda(Z_0) = \sum_T z^T, \tag{4.21}$$

where the sum is over all symplectic tableaux $T$ of shape $\lambda$.

Assume $q = \mathrm{Card}(Y)$ is sufficiently large, i.e. $q >> r$, and define functions $sc_\lambda(Z_0)$ by the identity

$$\sum_\lambda sc_\lambda(Z_0) s_\lambda(Y) = \frac{1}{\prod_{b \in B_0} \prod_{j=1}^q (1 - z_b y_j)} \prod_{1 \le i < j \le q} (1 - y_i y_j). \tag{4.22}$$

When $\lambda$ is a partition such that $\ell(\lambda) \le r$, then $sc_\lambda(Z_0) = sp_\lambda(Z_0)$ and these polynomials are the characters of the symplectic group $Sp(2r)$ and its Lie algebra $\mathfrak{sp}(2r)$. The combinatorics of these functions is discussed in [Su1].

Analogously suppose $Z_1 = \{z_b \mid b \in B_1\}$ and assume $q >> s$. Define functions $sb_\lambda(Z_1)$ by the identity

$$\sum_\lambda sb_\lambda(Z_1) s_\lambda(Y) = \frac{1}{\prod_{b \in B_1} \prod_{j=1}^q (1 - z_b y_j)} \prod_{1 \le i \le j \le q} (1 - y_i y_j). \tag{4.23}$$

When Card$(B_1)$ is odd and $\lambda$ is a partition such that $\ell(\lambda) \le s$, then the polynomials $sb_\lambda(Z_1)$ describe the characters of the orthogonal group $SO(2s+1)$ and its Lie algebra $\mathfrak{so}(2s+1)$. In this special case these functions have a combinatorial description similar to that of the symplectic Schur functions (see for example, [Su2]).

**Theorem 4.24.** *Let $\tilde{Z}$ be the set of variables $\{z_b\}$ indexed by the elements $b \in B = B_0 \cup B_1$ and let $Y = y_1, \ldots, y_q$ be an auxiliary set of variables such that $q >> r, s$. Then the following identities give equivalent definitions of polynomials $sc_\lambda(\tilde{Z})$:*

*(a)* $sc_\lambda(\tilde{Z}) = \sum_{\mu\subseteq\lambda} \left( \sum_{\pi} (-1)^{|\pi|/2} c^{\lambda}_{\mu,\pi} \right) s_\mu(\tilde{Z}),$

*where $s_\mu(\tilde{Z})$ denotes the hook Schur function of shape $\mu$, and $c^\lambda_{\mu,\pi}$ is the Littlewood-Richardson coefficient. The inner sum is over all partitions of the form $\pi = (r_1-1,\ldots,r_p-1 \mid r_1,\ldots,r_p)$, $r_1 \le q-1$, in Frobenius notation.*

*(b)* $sc_\lambda(\tilde{Z}) = \sum_{\tau\subseteq\lambda'} \left( \sum_{\rho} (-1)^{|\rho|/2} c^{\lambda'}_{\tau,\rho} \right) s_{\tau'}(\tilde{Z}),$

*where the inner sum is over all partitions of the form $\rho = (r_1+1,\ldots,r_p+1 \mid r_1,\ldots,r_p)$ , $r_1 \le q-1$, in Frobenius notation.*

*(c)* $\sum_\lambda sc_\lambda(\tilde{Z}) s_\lambda(Y) = \dfrac{\prod_{b\in B_1}\prod_{j=1}^q (1+z_b y_j)}{\prod_{b\in B_0}\prod_{j=1}^q (1-z_b y_j)} \prod_{1\le i<j\le q} (1-y_i y_j).$

*(d)* $\sum_\lambda sc_\lambda(\tilde{Z}) s_{\lambda'}(Y) = \dfrac{\prod_{b\in B_0}\prod_{j=1}^q (1+z_b y_j)}{\prod_{b\in B_1}\prod_{j=1}^q (1-z_b y_j)} \prod_{1\le i\le j\le q} (1-y_i y_j).$

*(e)* $sc_\lambda(\tilde{Z}) = \frac{1}{2}\det\left(h_{\lambda_i-i-j+2}(\tilde{Z}) + h_{\lambda_i-i+j}(\tilde{Z})\right).$

*where $h_\ell(\tilde{Z}) = s_{(\ell)}(\tilde{Z})$, and $(\ell)$ is the partition of $\ell$ with just one part.*

*(f)* $sc_\lambda(\tilde{Z}) = \frac{1}{2}\det\Big((e_{\lambda'_i-i-j+2}(\tilde{Z}) - e_{\lambda'_i-i-j}(\tilde{Z})) + (e_{\lambda'_i-i+j}(\tilde{Z}) - e_{\lambda'_i-i+j-2}(\tilde{Z}))\Big),$

*where $e_\ell(\tilde{Z}) = s_{(1^\ell)}(\tilde{Z})$, and $(1^\ell)$ is the partition of $\ell$ with all parts equal to 1.*

*(g)* $sc_\lambda(\tilde{Z}) = \det\left(e_{\lambda'_i-i+j}(\tilde{Z}) - e_{\lambda'_i-i-j}(\tilde{Z})\right).$

*(h)* $sc_\lambda(\tilde{Z}) = \sum_{\mu\subseteq\lambda} sc_\mu(Z_0) s_{\lambda'/\mu'}(Z_1),$

*where $s_{\lambda'/\mu'}(Z_1)$ is the skew Schur function in the variables $z_b$, $b \in B_1$.*

*(i)* $sc_\lambda(\tilde{Z}) = \sum_{\mu\subseteq\lambda} sb_{\mu'}(Z_1) s_{\lambda/\mu}(Z_0),$

*where $s_{\lambda/\mu}(Z_0)$ is the skew Schur function in the variables $z_b$, $b \in B_0$.*

*Proof.* Let us begin by proving (a) $\Longleftrightarrow$ (b), which follows from the calculation,

$$sc_\lambda(\tilde{Z}) = \sum_\mu \left( \sum_\pi (-1)^{|\pi|/2} c^\lambda_{\mu,\pi} \right) s_\mu(\tilde{Z}) = \sum_\mu \left( \sum_{\pi'} (-1)^{|\pi'|/2} c^{\lambda'}_{\mu',\pi'} \right) s_\mu(\tilde{Z})$$
$$= \sum_\mu \left( \sum_\rho (-1)^{|\rho|/2} c^{\lambda'}_{\mu',\rho} \right) s_\mu(\tilde{Z}) = \sum_\tau \left( \sum_\rho (-1)^{|\rho|/2} c^{\lambda'}_{\tau,\rho} \right) s_{\tau'}(\tilde{Z}).$$

Next we prove (c) using (a):

$$\begin{aligned}
\sum_{\lambda} sc_\lambda(\tilde{Z})s_\lambda(Y) &= \sum_{\lambda,\mu}\left(\sum_{\pi}(-1)^{|\pi|/2}c^{\lambda}_{\mu,\pi}\right)s_\mu(\tilde{Z})s_\lambda(Y)\\
&= \sum_{\mu} s_\mu(\tilde{Z})\sum_{\lambda,\pi}(-1)^{|\pi|/2}c^{\lambda}_{\mu,\pi}s_\lambda(Y)\\
&= \sum_{\mu} s_\mu(\tilde{Z})s_\mu(Y)\sum_{\pi}(-1)^{|\pi|/2}s_\pi(Y) && \text{by (4.9)},\\
&= \left(\sum_{\mu} s_\mu(\tilde{Z})s_\mu(Y)\right)\prod_{1\le i<j\le q}(1-y_iy_j) && \text{by (4.13)},\\
&= \frac{\prod_{b\in B_1}\prod_{j=1}^{q}(1+z_by_j)}{\prod_{b\in B_0}\prod_{j=1}^{q}(1-z_by_j)}\prod_{1\le i<j\le q}(1-y_iy_j) && \text{by (4.17)}.
\end{aligned}$$

The proof of (d) using (b) is virtually the same using (4.11) and (4.16).

Now let us argue that (c) implies (e): For each partition $\lambda = (\lambda_1,\dots,\lambda_q)$, $\lambda_1 \ge \lambda_2 \ge \cdots \ge \lambda_q \ge 0$ define

$$a_{\lambda+\delta}(Y) = \det(y_i^{\lambda_j+q-j}),$$

where the matrix $(y_i^{\lambda_j+q-j})$ is $q\times q$. Then $a_\delta(Y) = \prod_{1\le i<j\le q}(y_i-y_j)$, and the Schur function $s_\lambda(Y)$ can be written as

$$s_\lambda(Y) = \frac{a_{\lambda+\delta}(Y)}{a_\delta(Y)},$$

([M, Chap. I, §3.]). Let us write $y^{\lambda+\delta} = y_1^{\lambda_1+q-1}y_2^{\lambda_2+q-2}\cdots y_q^{\lambda_q}$. If $f(Y)$ is any symmetric function in the variables $y_1,\dots,y_q$, then the coefficient of $s_\lambda(Y)$ in the expansion of $f(Y)$ in terms of Schur functions is given by

$$[f(Y)]_{s_\lambda(Y)} = [f(Y)a_\delta(Y)]_{y^{\lambda+\delta}}.$$

Thus,

$$\begin{aligned}
sc_\lambda(\tilde{Z}) &= \left[\sum_{\lambda} sc_\lambda(\tilde{Z})s_\lambda(Y)\right]_{s_\lambda(Y)}\\
&= \left[\sum_{\lambda} sc_\lambda(\tilde{Z})s_\lambda(Y)a_\delta(Y)\right]_{y^{\lambda+\delta}}\\
&= \left[\frac{\prod_{b\in B_1}\prod_{j=1}^{q}(1+z_by_j)}{\prod_{b\in B_0}\prod_{j=1}^{q}(1-z_by_j)}\prod_{i<j}(1-y_iy_j)\prod_{i<j}(y_i-y_j)\right]_{y^{\lambda+\delta}},
\end{aligned}$$

where

$$\begin{aligned}\prod_{i<j}(1-y_iy_j)\prod_{i<j}(y_i-y_j) &= \prod_{i<j}(y_i-y_j)(1-y_iy_j)\\ &= \prod_{i<j}\Big(((y_j+y_j^{-1})-(y_i+y_i^{-1}))y_iy_j\Big)\\ &= (y_1\cdots y_q)^{q-1}\det\Big((y_i+y_i^{-1})^{q-(q-j+1)}\Big)\\ &= \tfrac{1}{2}(y_1\cdots y_q)^{q-1}\det\Big(y_i^{j-1}+y_i^{-(j-1)}\Big)\\ &= \tfrac{1}{2}\det\Big(y_i^{q+j-2}+y_i^{q-j}\Big).\end{aligned}$$

Now from (4.17) we see that

$$\frac{\prod_{b\in B_1}(1+z_by_i)}{\prod_{b\in B_0}(1-z_by_i)} = \sum_\lambda s_\lambda(\tilde{Z})s_\lambda(y_i) = \sum_{r\geq 0} s_{(r)}(\tilde{Z})s_{(r)}(y_i) = \sum_{r\geq 0} h_r(\tilde{Z})y_i^r.$$

Thus, we obtain

$$\begin{aligned} sc_\lambda(\tilde{Z}) &= \left[\tfrac{1}{2}\det\Big(y_i^{q+j-2}+y_i^{q-j}\Big)\prod_{i=1}^{q}\left(\sum_{r\geq 0}h_r(\tilde{Z})y_i^r\right)\right]_{y^{\lambda+\delta}}\\ &= \tfrac{1}{2}\det\left(\left[\sum_{r\geq 0}h_r(\tilde{Z})(y_i^{q+j-2+r}+y_i^{q-j+r})\right]_{y_i^{\lambda_i+q-i}}\right)\\ &= \tfrac{1}{2}\det\Big(h_{\lambda_i-i-j+2}(\tilde{Z})+h_{\lambda_i-i+j}(\tilde{Z})\Big).\end{aligned}$$

Identity (f) is gotten from (d) by a similar argument.

To show (f) implies (g), add columns $j-2, j-4, j-6, \ldots$ to column $j$ for each $j$. That is, add column 1 to column 3, add column 2 to column 4, add columns 1 and 3 to column 5, add columns 2 and 4 to column 6 and so on. The result follows immediately.

To derive (h) from (c), just compare the coefficients of $s_\lambda(Y)$ on both sides of the following equation:

$$\begin{aligned}\sum_\lambda sc_\lambda(\tilde{Z})s_\lambda(Y) &= \frac{\prod_{b\in B_1}\prod_{j=1}^q(1+z_by_j)}{\prod_{b\in B_0}\prod_{j=1}^q(1-z_by_j)}\prod_{i<j}(1-y_iy_j)\\ &= \left(\sum_\mu sc_\mu(Z_0)s_\mu(Y)\right)\left(\sum_\nu s_{\nu'}(Z_1)s_\nu(Y)\right) && \text{by (4.22)}\\ &= \sum_\lambda s_\lambda(Y)\sum_\mu sc_\mu(Z_0)\sum_\nu c^\lambda_{\mu,\nu}s_{\nu'}(Z_1) && \text{by (4.9)}\\ &= \sum_\lambda s_\lambda(Y)\sum_\mu sc_\mu(Z_0)\sum_\nu c^{\lambda'}_{\mu',\nu'}s_{\nu'}(Z_1) && \text{by (4.10)}\\ &= \sum_\lambda s_\lambda(Y)\sum_\mu sc_\mu(Z_0)s_{\lambda'/\mu'}(Z_1) && \text{by (4.11).}\end{aligned}$$

Obtaining (i) from (d) just amounts to a similar argument comparing the coefficients of $s_\lambda(Y)$ on both sides of the following equation:

$$\sum_\lambda sc_{\lambda'}(\tilde Z)s_\lambda(Y) = \frac{\prod_{b\in B_0}\prod_{j=1}^q(1+z_by_j)}{\prod_{b\in B_1}\prod_{j=1}^q(1-z_by_j)}\prod_{1\le i\le j\le q}(1-y_iy_j). \qquad \blacksquare$$

**Theorem 4.25.** *Assume $|n-m|>k$, and for $\lambda\vdash k-2h$, $h=0,1,\dots\lfloor k/2\rfloor$, let $T^\lambda$ be the $\mathfrak{spo}(V,\beta)$-module in (4.1). Let $sc_\lambda(\tilde Z)$ be the function defined in Theorem 4.24 and let $\mathrm{char}(T^\lambda)$ be as defined in Lemma 4.7 (d). Then*

$$\mathrm{char}(T^\lambda) = sc_\lambda(\tilde Z).$$

*Proof.* Suppose $\mu\vdash k-2j$ for some $j=0,1,\dots\lfloor k/2\rfloor$. Then by Lemma 4.7 (d) and Lemma 4.8,

$$(n-m)^j\overline{p}_\mu(\tilde Z) = \sum_{\lambda\in\widehat{B}_k}\chi^{\lambda'}_{B_k(n-m)}(e^{\otimes j}\otimes\gamma_\mu)\mathrm{char}(T^\lambda).$$

Thus, it follows from Theorem 4.4 that

$$\overline{p}_\mu(\tilde Z) = \sum_{\lambda\in\widehat{B}_k}\sum_{\pi\supseteq\lambda'}\left(\sum_{\theta\text{ even}} c^\pi_{\lambda',\theta}\right)\chi^\pi_{S_{k-2j}}(\mu)\mathrm{char}(T^\lambda).$$

By (4.20) and the orthogonality of characters for the symmetric group we have, for $\nu\vdash k-2j$,

$$\begin{aligned}
s_\nu(\tilde Z) &= \sum_{\mu\vdash k-2j}\sum_{\lambda\in\widehat{B}_k}\sum_{\pi\supseteq\lambda'}\left(\sum_{\theta\text{ even}} c^\pi_{\lambda',\theta}\right)\frac{\chi^{\nu'}_{S_{k-2j}}(\mu)\chi^\pi_{S_{k-2j}}(\mu)}{\mathfrak{z}_\mu}\mathrm{char}(T^\lambda)\\
&= \sum_{\lambda\in\widehat{B}_k}\sum_{\pi\supseteq\lambda'}\left(\sum_{\theta\text{ even}} c^\pi_{\lambda',\theta}\right)\mathrm{char}(T^\lambda)\delta_{\nu',\pi} = \sum_{\lambda\in\widehat{B}_k}\left(\sum_{\theta\text{ even}} c^{\nu'}_{\lambda',\theta}\right)\mathrm{char}(T^\lambda).
\end{aligned}$$

Thus, by Theorem 4.24 (b) and (4.16) it follows that

$$\begin{aligned}
sc_\tau(\tilde Z) &= \sum_{\nu\subseteq\tau'}\left(\sum_\rho(-1)^{|\rho|/2}c^{\tau'}_{\nu,\rho}\right)s_{\nu'}(\tilde Z)\\
&= \sum_{\lambda\in\widehat{B}_k}\sum_{\nu\subseteq\tau'}\left(\sum_\rho(-1)^{|\rho|/2}c^{\tau'}_{\nu,\rho}\right)\left(\sum_{\theta\text{ even}} c^{\nu}_{\lambda',\theta}\right)\mathrm{char}(T^\lambda)\\
&= \delta_{\tau',\lambda'}\mathrm{char}(T^\lambda),
\end{aligned}$$

for all partitions $\lambda\vdash k-2h$, $h=0,1,\dots,\lfloor k/2\rfloor$. $\blacksquare$

*Remark 4.26.* It follows immediately from the definition of the skew Schur functions, the definition of $sp_\lambda(Z_0)$ in (4.21), and Theorem 4.24(h) that the functions $sc_\lambda(\tilde Z)$ for certain partitions $\lambda$ can be expressed as a sum over weights of tableaux which have a symplectic part and a row-strict part. It is this interpretation which allows us to develop an insertion scheme modelling these functions in the next section.

## 5. TABLEAUX AND AN INSERTION SCHEME

*Up-down tableaux and $\mathfrak{spo}$-tableaux*

Fix positive integers $r$ and $n$. An *up-down tableau of length $k$ and shape $\lambda$* is a sequence of partitions $\Lambda = (\lambda^0, \lambda^1, \ldots, \lambda^k)$ such that $\lambda^0 = \emptyset$ and $\lambda^k = \lambda$ and, for each $i = 1, \ldots, k$, $\lambda^i$ is obtained from $\lambda^{i-1}$ by either adding or removing a box. An *up-down $(r,n)$-tableau* is an up-down tableau $\Lambda = (\emptyset = \lambda^0, \lambda^1, \ldots, \lambda^k = \lambda)$ of length $k$ and shape $\lambda$ such that each $\lambda^i$ is an *$(r,n)$-hook shape partition*, that is, $\lambda^i_{r+1} \leq n$. If $\lambda^k = \lambda$ is the final shape, then $\lambda \vdash k - 2h$ for some $0 \leq h \leq \lfloor k/2 \rfloor$.

**Example.** When $k = 7$,

$\Lambda =$ 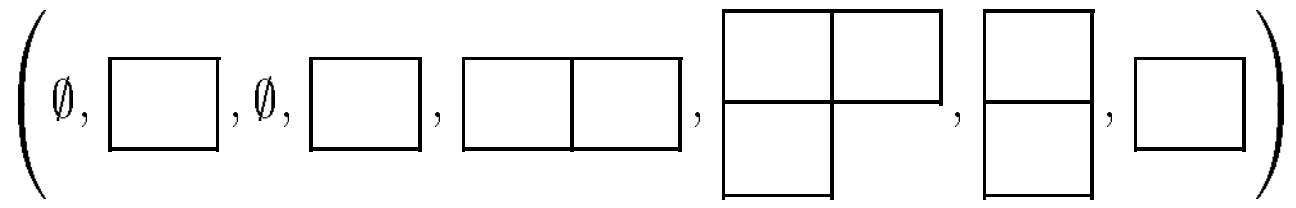

is an up-down $(r,n)$-tableau of shape $\lambda = (1, 0, \ldots)$ for all $r \geq 1$ and $n \geq 1$.

Let $B_0 = \{t_1, t_1^*, \ldots, t_r, t_r^*\}$ and $B_1 = \{v_1, \ldots, v_n\}$ and let $B = B_0 \cup B_1$. Here we do not need to distinguish between $n$ odd or $n$ even and for that reason we do not use the notation $u_1, u_1^*, \ldots, u_s, u_s^*, (u_{s+1})$ for the elements of $B_1$ as we have done in previous sections. Order $B$ as follows:

$$B = \{t_1 < t_1^* < t_2 < t_2^* < \ldots < t_r < t_r^* < v_1 < \ldots < v_n\}.$$

Suppose $m = 2r$ and $n$ are fixed. An *$\mathfrak{spo}(m,n)$-tableau $T$ of shape $\lambda$* (or simply an $\mathfrak{spo}$-tableau for short) is a filling of the boxes in the Ferrers diagram of $\lambda$ with entries from $B$ such that

($\mathfrak{spo}$.1) the subtableau $S$ of $T$ obtained by taking all the boxes with entries from $B_0$ is a column-strict tableau of partition shape (its entries are weakly increasing from left to right across each row and are strictly increasing from top to bottom down each column), and the entries in row $i$ are $\geq t_i$ for each row in $S$;

($\mathfrak{spo}$.2) the skew tableau $T/S$ is row-strict (its entries are strictly increasing from left to right across each row and are weakly increasing from top to bottom down each column).

In an $\mathfrak{spo}(m,n)$-tableau the entries in row $i$ for $i \geq r+1$ necessarily belong to $B_1$, and since the skew tableau $T/S$ is row-strict, the underlying partition must be an $(r,n)$-hook shape partition.

**Example.** Let

$$T = \begin{array}{|c|c|c|c|c|}\hline t_1 & t_1^* & t_2 & v_3 & v_4 \\ \hline t_2^* & t_2^* & v_1 & v_3 & \multicolumn{1}{c}{} \\ \cline{1-4} t_3 & t_3 & \multicolumn{3}{c}{} \\ \cline{1-2} t_4 & \multicolumn{4}{c}{} \\ \cline{1-1} \end{array} \, .$$

Then $T$ is an $\mathfrak{spo}(m,n)$-tableau for all $m \geq 8$ and $n \geq 4$ where

$$S = \begin{array}{|c|c|c|}\hline t_1 & t_1^* & t_2 \\ \hline t_2^* & t_2^* & \multicolumn{1}{c}{} \\ \cline{1-2} t_3 & t_3 & \multicolumn{1}{c}{} \\ \cline{1-2} t_4 & \multicolumn{2}{c}{} \\ \cline{1-1} \end{array} \quad \text{and} \quad T/S = \begin{array}{|c|c|c|c|c|}\hline \phantom{t} & \phantom{t} & \phantom{t} & v_3 & v_4 \\ \hline & & v_1 & v_3 & \multicolumn{1}{c}{} \\ \cline{1-4} & & \multicolumn{3}{c}{} \\ \cline{1-2} & \multicolumn{4}{c}{} \\ \cline{1-1} \end{array} \, .$$

**Theorem 5.1.** *If $T$ is an $\mathfrak{spo}(m,n)$-tableau, replace each entry by the corresponding variable $z_b$ as defined in (4.5) and let $z^T$ be the product of the entries in the result. Let $sc_\lambda(\tilde{Z})$ be the polynomial defined in Theorem 4.24. Then*

$$sc_\lambda(\tilde{Z}) = \sum_T z^T$$

*where the sum is over all the $\mathfrak{spo}(m,n)$-tableaux of shape $\lambda$.*

*Proof.* Suppose $T$ is an $\mathfrak{spo}(m,n)$-tableau, and let $S$ denote the subtableau of $T$ containing the entries $t_1, t_1^*, \ldots, t_r, t_r^*$. Necessarily the number of rows in $S$ is $\leq r$. When the entries $t_i, t_i^*$, $i = 1, \ldots, r$ of $S$ are replaced by the elements $z_{t_i}, z_{t_i^*}$, then the result is a symplectic tableau as defined just before equation (4.21). Likewise, if we substitute the element $z_{v_j}$ for the entry $v_j$, $j = 1, \ldots, n$, in $T/S$ and transpose the resulting tableau, what is produced is a column-strict tableau. The result then follows from the definition of the skew Schur function, (4.21), and Theorem 4.24(h). ∎

Assume

$\mathcal{W}_k$ is the set of words of length $k$ in the alphabet $B$, and

$\mathcal{P}_k$ is the set of pairs $(T, \Lambda)$ consisting of a $\mathfrak{spo}(m,n)$-tableau $T$ of shape $\lambda$ and an up-down $(r,n)$-tableau $\Lambda = (\lambda^0, \lambda^1, \ldots, \lambda^k)$ of length $k$ and shape $\lambda$.

We shall prove in Theorem 5.5 that there is a bijection between $\mathcal{W}_k$ and $\mathcal{P}_k$. Accomplishing this requires introducing a few more definitions.

*Punctured tableaux and the maps jeu and injeu*

A *punctured tableau* is a tableau with exactly one empty box. A *punctured $\mathfrak{spo}$-tableau* is a punctured tableau obtained from an $\mathfrak{spo}$-tableau by removing the entry from exactly one box. If a partition $\lambda$ contains a box at position $(i,j)$, and there is no box at locations $(i, j+1)$ and $(i+1, j)$ in $\lambda$, then $(i, j)$ is said to be a *corner* of $\lambda$. When a punctured tableau has its empty box at a corner, we may identify the punctured tableau with the tableau obtained by removing the empty box.

In Lemma 5.2 below we shall define operations $se$ and $nw$ on the set of punctured $\mathfrak{spo}$-tableaux, but first we make the following definitions.

(i) $\mathcal{SE}$ is the set of punctured $\mathfrak{spo}$-tableaux whose empty box occurs at a corner.

(ii) $\mathcal{NW}$ is the set of punctured $\mathfrak{spo}$-tableaux whose empty box is in the first column, and either the empty box is in the first row or the tableau obtained by switching the empty box with its neighbor to the north (the one immediately above it) is not a punctured $\mathfrak{spo}$-tableau.

**Lemma 5.2.** *(Compare [Be, Lemma 2].) Let $T$ be a punctured $\mathfrak{spo}$-tableau.*

*(i) Suppose $T \notin \mathcal{SE}$. Then there is a unique way of switching the empty box of $T$ either with its neighbor to the east (the box on its right) or with its neighbor to the south (the box immediately below) so that the resulting tableau is a punctured $\mathfrak{spo}$-tableau. We denote the resulting tableau by $se(T)$.*

*(ii) Suppose $T \notin \mathcal{NW}$. Then there is a unique way of switching the empty box of $T$ either with its neighbor to the north or with its west neighbor (the box on its left) so that the resulting tableau is a punctured $\mathfrak{spo}$-tableau. We denote this resulting tableau by $nw(T)$.*

*(iii) If $T \notin \mathcal{SE}$, then $nw(se(T)) = T$, and if $T \notin \mathcal{NW}$, then $se(nw(T)) = T$.*

*Proof.* Assume the neighbors of the empty box are as pictured below, where one or more of the boxes with letters may be absent. (We are not necessarily displaying the entire tableau rather just the neighbors of the empty box since they are all that matter in the argument.)

$$T = \begin{array}{|c|c|c|} \hline x & a & y \\ \hline b & & c \\ \hline z & d & w \\ \hline \end{array}$$

First consider the case that all the neighbors of the empty box are present. We note that $a < c$ since $a \le y \le c$, and at least one of the inequalities is strict. Similarly $b < d$ holds since $b \le z \le d$ and at least one of the inequalities is strict.

To prove (i) consider the tableaux

$$e(T) = \begin{array}{|c|c|c|} \hline x & a & y \\ \hline b & c & \\ \hline z & d & w \\ \hline \end{array} \quad \text{and} \quad s(T) = \begin{array}{|c|c|c|} \hline x & a & y \\ \hline b & d & c \\ \hline z & & w \\ \hline \end{array}.$$

If $c = d \in B_0$, then only $s(T)$ can be an $\mathfrak{spo}$-tableau. If $d$ is in row $i$ of $T$, then since $T$ is an $\mathfrak{spo}$-tableau, $d \ge t_i > t_{i-1}$ so that condition ($\mathfrak{spo}$.1) is fulfilled in $s(T)$. If $c = d \in B_1$, then only $e(T)$ is an $\mathfrak{spo}$-tableau. If $c > d$, then only $s(T)$ is an $\mathfrak{spo}$-tableau, while if $c < d$, then just $e(T)$ is an $\mathfrak{spo}$-tableau. The tableau $se(T)$ is $s(T)$ or $e(T)$, whichever one is an $\mathfrak{spo}$-tableau. (Note $se(T)$ should not be confused with $s(e(T))$.)

For (ii) consider the tableaux

$$w(T) = \begin{array}{|c|c|c|}\hline x & a & y \\ \hline & b & c \\ \hline z & d & w \\ \hline \end{array} \quad \text{and} \quad n(T) = \begin{array}{|c|c|c|}\hline x & & y \\ \hline b & a & c \\ \hline z & d & w \\ \hline \end{array}.$$

If $b = a \in B_0$ and $b$ is in row $i$, then $t_i \leq b = a < d$ and only $n(T)$ is an $\mathfrak{spo}$-tableau. If $b = a \in B_1$, then $b = a < c$, and $w(T)$ is an $\mathfrak{spo}$-tableau but $n(T)$ is not. When $a > b$, then just $n(T)$ is an $\mathfrak{spo}$-tableau, and when $a < b$, then only $w(T)$ is an $\mathfrak{spo}$-tableau. The tableau $nw(T)$ is either $n(T)$ or $w(T)$, whichever happens to be an $\mathfrak{spo}$-tableau.

To prove (iii) suppose that $T \notin \mathcal{SE}$ and $se(T) = s(T)$. By (ii) $nw(s(T))$ is either $n(s(T))$ or $w(s(T))$ depending on which is an $\mathfrak{spo}$-tableau, and clearly $n(s(T)) = T$ is an $\mathfrak{spo}$-tableau. Thus, $nw(se(T)) = n(s(T)) = T$ in this case. Similarly, if $se(T) = e(T)$, then $nw(s(T)) = w(e(T)) = T$. For $T \notin \mathcal{NW}$, using (i) we see that $se(nw(T)) = T$.

Finally, we note that when the empty box of $T$ has only one of its neighbors to the east or south, arguments similar to ones above can be used to show that $se(T)$ is obtained by switching the empty box with the available neighbor. Also, when the empty box of $T$ has just one of its northern or western neighbors $nw(T)$ is obtained by switching the empty box with the available neighbor. In these cases, (iii) holds similarly. ∎

Let $\mathcal{SE}$ and $\mathcal{NW}$ be the sets of punctured $\mathfrak{spo}$-tableaux described just before Lemma 5.2, and define maps $jeu: \mathcal{NW} \to \mathcal{SE}$ and $injeu: \mathcal{SE} \to \mathcal{NW}$ as follows:

(i) For $T \in \mathcal{NW}$, there is a least integer $j \geq 1$ such that applying $se$ to $T$ $j$ times gives a tableau in $\mathcal{SE}$. Then $jeu(T) \stackrel{\text{def}}{=} (se)^j(T)$.

(ii) For $T \in \mathcal{SE}$, there is a least integer $i \geq 1$ such that $(nw)^i(T) \in \mathcal{NW}$. Then $injeu(T) \stackrel{\text{def}}{=} (nw)^i(T)$.

The following is an immediate consequence of Lemma 5.2.

**Lemma 5.3.** *The map jeu is a bijection from $\mathcal{NW}$ to $\mathcal{SE}$ with inverse mapping injeu.*

**Example.** The sequence below illustrates applying the mapping $jeu$ to a tableau in $\mathcal{NW}$ and ending with a tableau in $\mathcal{SE}$.

$$\begin{array}{|c|c|c|}\hline t_1^* & t_2 & t_2 \\ \hline & t_3 & v_2 \\ \hline t_3 & v_1 & v_2 \\ \hline v_2 \\ \cline{1-1} \end{array} \xrightarrow[se]{} \begin{array}{|c|c|c|}\hline t_1^* & t_2 & t_2 \\ \hline t_3 & t_3 & v_2 \\ \hline & v_1 & v_2 \\ \hline v_2 \\ \cline{1-1} \end{array} \xrightarrow[se]{} \begin{array}{|c|c|c|}\hline t_1^* & t_2 & t_2 \\ \hline t_3 & t_3 & v_2 \\ \hline v_1 & & v_2 \\ \hline v_2 \\ \cline{1-1} \end{array} \xrightarrow[se]{} \begin{array}{|c|c|c|}\hline t_1^* & t_2 & t_2 \\ \hline t_3 & t_3 & v_2 \\ \hline v_1 & v_2 & \\ \hline v_2 \\ \cline{1-1} \end{array}.$$

If we take the rightmost tableau in this sequence as the initial tableau and perform $nw$ successively, then the sequence moving from right to left corresponds to applying $injeu$.

*Insertion of a letter into an $\mathfrak{spo}$-tableau*

Let $T$ be an $\mathfrak{spo}$-tableau, and assume $a \in B$. We define an algorithm consisting of a sequence of steps which inserts $a$ into $T$ to yield a tableau $(a \rightarrow T)$.

(1) Start with $b = a$ and $i = j = 1$.

(2) If $b \in B_0$, then insert $b$ into the $i$th row of $T$ as follows: If there is an entry in row $i$ which is greater than $b$, then displace the leftmost such entry and insert $b$ into its box except in the following case. If $b = t_i$ and there is an $t_i^*$ in the $i$th row, then replace the leftmost $t_i^*$ in the row with $t_i$ and remove the entry in the $(i,1)$-position (which is necessarily an $t_i$) making it an empty box. If there is no entry in the $i$th row which is greater than $b$, then adjoin $b$ to the end of the row.

If $b \in B_1$, then insert $b$ into the $j$th column as follows: If there is an entry in the $j$th column which is greater than $b$, then displace the topmost such entry and insert $b$ into its position. If there is no entry in the $j$th column which is greater than $b$, then adjoin $b$ at the end of the column.

(3) Set $b$ equal to the displaced entry and change $i$ to $p+1$ and $j$ to $q+1$ where $(p,q)$ was the position of the displaced entry. Repeat step (2) until an entry is adjoined to the end of a row or a column, or an empty box is created.

(4) Let $(a \rightarrow T)'$ be the result of steps (1)-(3). Set $(a \rightarrow T) = (a \rightarrow T)'$ if $(a \rightarrow T)'$ is a tableau, and $(a \rightarrow T) = jeu((a \rightarrow T)')$ if $(a \rightarrow T)'$ is a punctured tableau.

**Lemma 5.4.** *Let $T$ be an $\mathfrak{spo}$-tableau and assume $a \in B$. Then $(a \rightarrow T)$ is an $\mathfrak{spo}$-tableau.*

*Proof.* The algorithm for insertion of a letter into an $\mathfrak{spo}$-tableau is just Berele's insertion algorithm (see [Be, Lemma 4]) for the part of the tableau involving entries from $B_0$, and it amounts to the usual Robinson-Schensted algorithm (see [S, §3.3], for example) for the skew-tableau with entries from $B_1$. Consequently, the result $(a \rightarrow T)'$ is either a tableau or a punctured tableau in $\mathcal{NW}$. From Lemma 5.2 it is clear that $(a \rightarrow T)$ is an $\mathfrak{spo}$-tableau provided we identify a punctured tableau in $\mathcal{SE}$ with the tableau obtained by removing its corner empty box. ∎

*Insertion of a word into an $\mathfrak{spo}$-tableau*

Recall that

$\mathcal{W}_k$ is the set of words of length $k$ in the alphabet $B$ and

$\mathcal{P}_k$ is the set of pairs $(T,\Lambda)$ consisting of a $\mathfrak{spo}(m,n)$-tableau $T$ of shape $\lambda$ and an up-down $(r,n)$-tableau $\Lambda = (\lambda^0, \lambda^1, \ldots, \lambda^k)$ of length $k$ and shape $\lambda$.

For $k \geq 0$, define maps $\mathfrak{spo}_k: \mathcal{W}_k \rightarrow \mathcal{P}_k$ inductively by

(i) $\mathfrak{spo}_0(w) = (\emptyset, (\emptyset))$, where $w$ is the emptyword;

(ii) If $w = w_1 \cdots w_k$ is a word of length $k$ and $\mathfrak{spo}_{k-1}(w_1 \cdots w_{k-1}) = (T^{k-1}, (\lambda^0, \lambda^1, \ldots, \lambda^{k-1}))$, then $\mathfrak{spo}_k(w) = (T^k, (\lambda^0, \lambda^1, \ldots, \lambda^k))$ where $T^k = (w_k \rightarrow T^{k-1})$ and $\lambda^k$ is the underlying partition of $T^k$.

**Theorem 5.5.** *The map $\mathfrak{spo}_k : \mathcal{W}_k \longrightarrow \mathcal{P}_k$ is a bijection.*

*Proof.* It is clear that $\mathfrak{spo}_k(w) \in \mathcal{P}_k$ for each $w \in \mathcal{W}_k$. To show that $\mathfrak{spo}_k$ is a bijection, it is enough to find its inverse, and for that it suffices to find the inverse of the step $(w_\ell \rightarrow$ $T^{\ell-1}) = T^\ell$. Let $\lambda^{\ell-1}$ and $\lambda^\ell$ be the underlying partitions of $T^{\ell-1}$ and $T^\ell$, respectively. We present an algorithm to produce $w_\ell$ from $\lambda^{\ell-1}$, $\lambda^\ell$ and $T^\ell$:

(1) (initial step) If $\lambda^\ell$ is obtained from $\lambda^{\ell-1}$ by adding a box, and if the entry in the corresponding box of $T^\ell$ is $a$, then delete the adjoined box from $T^\ell$ and apply (2) below beginning with $b = a$.

If $\lambda^\ell$ is obtained from $\lambda^{\ell-1}$ by removing a box, then consider the punctured $\mathfrak{spo}$-tableau in $\mathcal{SE}$ obtained by adding an empty box to $T^\ell$ in the corresponding corner. Apply *injeu* to obtain a punctured $\mathfrak{spo}$-tableau in $\mathcal{NW}$. The empty box is now in position $(i,1)$ for some $i$. Fill the empty box with an $t_i$, and then change the last $t_i$ in row $i$ to an $t_i^*$. Apply (2) below beginning with $b = t_i$ regarded as displaced from position $(i,1)$.

(2) (general step) Suppose that an entry $b$ is displaced from the $(i,j)$-position of $T^\ell$.

If $b \in B_0$, then from the $(i-1)$st row of $T^\ell$ remove the rightmost entry which is smaller than $b$ and put $b$ in the empty box.

If $b \in B_1$, then from the $(j-1)$st column of $T^\ell$ remove the bottommost entry which is smaller than $b$ and put $b$ in the empty box.

Repeat this step until an entry in $B_0$ is displaced from the first row or an entry in $B_1$ is displaced from the first column. That entry is $w_\ell$.

This gives the desired inverse process, and hence completes the proof. ∎

We refer to the map $\mathfrak{spo}_k \colon \mathcal{W}_k \longrightarrow \mathcal{P}_k$ as $\mathfrak{spo}$-*insertion* and the inverse described in the proof of Theorem 5.5 as $\mathfrak{spo}$-*deletion*. We illustrate $\mathfrak{spo}$-insertion and $\mathfrak{spo}$-deletion in the following example.

**Example.** Suppose $r = 4$, $n = 4$, and $k = 6$. Then $B = \{t_1 < t_1^* < t_2 < t_2^* < t_3 < t_4 < t_4^* < v_1 < v_2 < v_3 < v_4\}$. Consider the word $w = v_2\ t_2\ t_1^*\ t_1\ v_3\ t_1$ of length 6. Then $\mathfrak{spo}$-insertion for the word $w$ is given by

$$T^0 = \emptyset, \quad T^1 = \begin{array}{|c|}\hline v_2 \\ \hline\end{array}, \quad T^2 = \begin{array}{|c|c|}\hline t_2 & v_2 \\ \hline\end{array}, \quad T^3 = \begin{array}{|c|c|}\hline t_1^* & v_2 \\ \hline t_2 \\ \cline{1-1}\end{array},$$

$$\begin{array}{|c|c|}\hline & v_2 \\ \hline t_2 \\ \cline{1-1}\end{array} \xrightarrow{jeu} \begin{array}{|c|c|}\hline t_2 & v_2 \\ \hline\end{array} = T^4, \quad T^5 = \begin{array}{|c|c|}\hline t_2 & v_2 \\ \hline v_3 \\ \cline{1-1}\end{array}, \quad T^6 = \begin{array}{|c|c|}\hline t_1 & v_2 \\ \hline t_2 & v_3 \\ \hline\end{array}.$$

Thus $\mathfrak{spo}(w) = (T, \Lambda)$ where $T = T^6$ and

$$\Lambda = \left( \emptyset, \begin{array}{|c|}\hline \ \\ \hline\end{array}, \begin{array}{|c|c|}\hline \ & \ \\ \hline\end{array}, \begin{array}{|c|c|}\hline \ & \ \\ \hline \ \\ \cline{1-1}\end{array}, \begin{array}{|c|c|}\hline \ & \ \\ \hline\end{array}, \begin{array}{|c|c|}\hline \ & \ \\ \hline \ \\ \cline{1-1}\end{array}, \begin{array}{|c|c|}\hline \ & \ \\ \hline \ & \ \\ \hline\end{array} \right).$$

To recover $w$ from $(T, \Lambda)$ we apply $\mathfrak{spo}$-deletion:

$$T_6 = \begin{array}{|c|c|}\hline t_1 & u_2 \\ \hline t_2 & v_3 \\ \hline\end{array} \longrightarrow \begin{array}{|c|c|}\hline t_1 & v_2 \\ \hline v_3 \\ \cline{1-1}\end{array} \longrightarrow \begin{array}{|c|c|}\hline t_2 & v_2 \\ \hline v_3 \\ \cline{1-1}\end{array} \quad \text{and} \quad w_6 = t_1,$$

$$T_5 = \begin{array}{|c|c|}\hline t_2 & v_2 \\ \hline v_3 & \multicolumn{1}{c}{} \\ \cline{1-1}\end{array} \longrightarrow \begin{array}{|c|c|}\hline t_2 & v_2 \\ \hline\end{array} \text{ and } w_5 = v_3,$$

$$T_4 = \begin{array}{|c|c|}\hline t_2 & v_2 \\ \hline\end{array} \longrightarrow \begin{array}{|c|c|}\hline t_2 & v_2 \\ \hline & \multicolumn{1}{c}{} \\ \cline{1-1}\end{array} \xrightarrow[inj\,eu]{} \begin{array}{|c|c|}\hline & v_2 \\ \hline t_2 & \multicolumn{1}{c}{} \\ \cline{1-1}\end{array} \longrightarrow \begin{array}{|c|c|}\hline t_1 & v_2 \\ \hline t_2 & \multicolumn{1}{c}{} \\ \cline{1-1}\end{array} \longrightarrow \begin{array}{|c|c|}\hline t_1^* & v_2 \\ \hline t_2 & \multicolumn{1}{c}{} \\ \cline{1-1}\end{array} \text{ and } w_4 = t_1,$$

$$T_3 = \begin{array}{|c|c|}\hline t_1^* & v_2 \\ \hline t_2 & \multicolumn{1}{c}{} \\ \cline{1-1}\end{array} \longrightarrow \begin{array}{|c|c|}\hline t_2 & v_2 \\ \hline\end{array} \text{ and } w_3 = t_1^*,$$

$$T_2 = \begin{array}{|c|c|}\hline t_2 & v_2 \\ \hline\end{array} \longrightarrow \begin{array}{|c|}\hline v_2 \\ \hline\end{array} \text{ and } w_2 = t_2,$$

$$T_1 = \begin{array}{|c|}\hline v_2 \\ \hline\end{array} \longrightarrow \emptyset \text{ and } w_1 = v_2.$$

Thus, we get $w = v_2\, t_2\, t_1^*\, t_1\, v_3\, t_1$.

**Corollary 5.6.** *For positive integers $m = 2r$ and $n$,*

$$(m+n)^k = \sum_{h=0}^{\lfloor k/2 \rfloor} \sum_{\lambda \vdash k-2h} \mathrm{ud}_\lambda(r,n) \cdot \mathrm{spo}_\lambda(m,n)$$

*where $\mathrm{ud}_\lambda(r,n)$ is the number of up-down $(r,n)$-tableaux of shape $\lambda$, and $\mathrm{spo}_\lambda(m,n)$ is the number of $\mathfrak{spo}(m,n)$-tableaux of shape $\lambda$.*

*Proof.* This is an immediate consequence of Theorem 5.5 by counting the number of elements in the sets $\mathcal{W}_k$ and $\mathcal{P}_k$. ∎

*Remark 5.7.* By ([We],[W]), the dimension of the irreducible module $B_\lambda$, $\lambda \in \widehat{B}_k = \{\lambda \vdash k-2h \mid h = 0,1,\ldots,\lfloor k/2 \rfloor\}$, for the Brauer algebra $B_k(n-m)$ is the number of up-down tableaux of length $k$ and shape $\lambda$ whenever $|n-m| > k$. This is also the number of up-down tableaux of length $k$ and shape $\lambda'$ whenever $|n-m| > k$. By Theorem 4.24, $\mathrm{char}(T^\lambda) = sc_\lambda(\tilde{Z})$. Thus, Theorem 5.5 gives an explicit bijection realizing the identity

$$V^{\otimes k} \cong \bigoplus_\lambda T^\lambda \otimes B_{\lambda'}$$

from Proposition 4.2 (b).